\documentclass[10pt,oneside,final]{amsart}
\usepackage[T1]{fontenc}
\usepackage[latin9]{inputenc}
\usepackage{geometry}
\usepackage{amssymb}
\usepackage{esint}

\makeatletter
  \theoremstyle{plain}
  \newtheorem{thm}{Theorem}
  \theoremstyle{remark}
  \newtheorem*{acknowledgement*}{Acknowledgement}
  \theoremstyle{plain}
  \newtheorem{lem}{Lemma}
  \theoremstyle{plain}
  \newtheorem{prop}{Proposition}
  \theoremstyle{remark}
  \newtheorem{rem}{Remark}
  \theoremstyle{plain}
  \newtheorem{cor}{Corollary}
  \theoremstyle{plain}
  \newtheorem{conjecture}{Conjecture}
 \theoremstyle{definition}
  \newtheorem{example}{Example}

\makeatother

\begin{document}

\title{Lower estimates on microstates free entropy dimension.}

\author{Dimitri Shlyakhtenko}

\address{Department of Mathematics, UCLA, Los Angeles, CA 90095.}

\email{shlyakht@math.ucla.edu}

\thanks{Research supported by NSF grant DMS-0555680}

\begin{abstract}
By proving that certain free stochastic differential equations with
analytic coefficients have stationary solutions, we give a lower estimate
on the microstates free entropy dimension of certain $n$-tuples $X_{1},\ldots,X_{n}$.
In particular, we show that $\delta_{0}(X_{1},\ldots,X_{n})\geq\dim_{M\bar{\otimes}M^{o}}V$
where $M=W^{*}(X_{1},\dots,X_{n})$ and $V=\{(\partial(X_{1}),\dots,\partial(X_{n})):\partial\in\mathcal{C}\}$
is the set of values of derivations $A=\mathbb{C}[X_{1},\ldots X_{n}]\to A\otimes A$
with the property that $\partial^{*}\partial(A)\subset A$. We show
that for $q$ sufficiently small (depending on $n$) and $X_{1},\ldots,X_{n}$
a $q$-semicircular family, $\delta_{0}(X_{1},\ldots,X_{n})>1$. In
particular, for small $q$, $q$-deformed free group factors have
no Cartan subalgebras. An essential tool in our analysis is a free
analog of an inequality between Wasserstein distance and Fisher information
introduced by Otto and Villani (and also studied in the free case
by Biane and Voiculescu).
\end{abstract}
\maketitle

\section{Introduction.}

We present in this paper a general technique for proving lower estimates
for Voiculescu's microstates free entropy dimension $\delta_{0}$.
The free entropy dimension $\delta_{0}$ was introduced in \cite{dvv:entropy2,dvv:entropy3}
and is a number associated to an $n$-tuple of self-adjoint elements
$X_{1},\ldots,X_{n}$ in a tracial von Neumann algebra. This quantity
was used by Voiculescu and others (see e.g \cite{dvv:entropy3,ge:entropy2,ge-shen:commutingGenerators,stephan:thinness,jung:onebounded})
to prove a number of very important results in von Neumann algebras.
These results often take the form: if $\delta_{0}(X_{1},\ldots,X_{n})>1$,
then $M=W^{*}(X_{1},\ldots,X_{n})$ cannot have certain decomposition
properties (e.g., is non-$\Gamma$, has no Cartan subalgebras, is
not a non-trivial tensor product and so on). For this reason, it is
important to know if some given von Neumann algebra has a set of generators
with the property that $\delta_{0}>1$. We prove that this is the
case (for small values of $q$) for the {}``$q$-deformed free group
factors'' of Bozejko and Speicher \cite{Speicher:q-comm}:

\begin{thm}
For a fixed $N$, and all $|q|<(4N^{3}+2)^{-1}$, the $q$-semicircular
family $X_{1},\ldots,X_{N}$ satisfies $\delta_{0}(X_{1},\ldots,X_{N})>1$
and $\delta_{0}(X_{1},\ldots,X_{N})\geq N(1-q^{2}N(1-q^{2}N)^{-1})$.
\end{thm}
The theorem applies for $|q|\leq0.029$ if $N=2$. Combined with the
available results on free entropy dimension, we obtain that in this
range of values of $q$, the algebras $\Gamma_{q}(\mathbb{R}^{N})=W^{*}(X_{1},\ldots,X_{N})$
have no Cartan subalgebras (or, more generally, that $\Gamma_{q}(\mathbb{R}^{N})$,
when viewed as a bimodule over any of its abelian subalgebras, contain
a coarse sub-bimodule). One also gets that these algebras are prime
(although this was already proved using Ozawa's techniques from \cite{ozawa:solid}
elsewhere \cite{shlyakht:qdim}).

The free entropy dimension $\delta_{0}$ is closely related to $L^{2}$
Betti numbers (see \cite{connes-shlyakht:l2betti,shlyakht-mineyev:freedim}),
more precisely, with Murray-von Neumann dimensions of spaces of certain
derivations. For example, the non-microstates free entropy dimension
$\delta^{*}$ (which is the non-microstates {}``relative'' of $\delta_{0}$)
is in many cases equal to $L^{2}$ Betti numbers of the underlying
(non-closed) algebra \cite{shlyakht-mineyev:freedim,shlyakht:counterexample}.
It is known that $\delta_{0}\leq\delta^{*}$ and thus it is important
to find lower estimates for $\delta_{0}$ in terms of dimensions of
spaces of derivations. To this end we prove:

\begin{thm}
\label{thm:estAlg}Let $(A,\tau)$ be a finitely-generated algebra
with a positive trace $\tau$ and generators $X_{1},\ldots,X_{N}$,
and let $\operatorname{Der}_{c}(A;A\otimes A)$ denote the space of
derivations from $A$ to $A\otimes A$ which are $L^{2}$ closable
and so that $\partial^{*}\partial(X_{j})\in A$. Consider the A,A-bimodule\[
V=\{(\delta(X_{1}),\ldots,\delta(X_{n})):\delta\in\operatorname{Der}_{c}(A;A\otimes A)\}\subset(A\otimes A)^{N}.\]
Assume finally that $M=W^{*}(A,\tau)$ can be embedded into the ultrapower
of the hyperfinite II$_{1}$ factor. Then\[
\delta_{0}(X_{1},\ldots,X_{n})\geq\dim_{M\bar{\otimes}M^{o}}\overline{V}^{L^{2}(A\otimes A,\tau\otimes\tau)^{N}}.\]

\end{thm}
We actually prove Theorem \ref{thm:estAlg} under a less restrictive
assumption: we require that $\delta(X_{j})$ and $\delta^{*}\delta(X_{j})$
be {}``analytic'' as functions of $X_{1},\dots,X_{N}$; more precisely,
there should exist non-commutative power series $\Xi_{j}$ and $\xi_{j}$
with sufficiently large multi-radii of convergence so that $\delta(X_{j})=\Xi_{j}(X_{1},\dots,X_{N})$
and $\delta^{*}\delta(X_{j})=\xi_{j}(X_{1},\dots,X_{N})$; see Theorem
\ref{thm:deltaEstimate} below for a precise statement. 

This theorem is a rich source of lower estimates for $\delta_{0}$.
For example, if $T\in A\otimes A$, then $\delta:X\mapsto[X,T]=XT-TX$
is a derivation in $\operatorname{Der}_{c}(A;A\otimes A)$. If $W^{*}(A)$
is diffuse, then the map\[
L^{2}(A\otimes A)\ni T\mapsto([T,X_{1}],\ldots,[T,X_{N}])\to L^{2}(A\otimes A)^{N}\]
is injective and thus the dimension over $M\bar{\otimes}M^{o}$ of
its image is the same as the dimension of $L^{2}(A\otimes A)$, i.e.,
$1$. Hence $\dim_{M\bar{\otimes}M^{o}}\overline{V}\geq1$ and so
$\delta_{0}(X_{1},\ldots,X_{n})\geq1$ if $W^{*}(A)$ is $R^{\omega}$
embeddable ({}``hyperfinite monotonicity'' of \cite{jung-freexentropy}). 

If the two tuples $X_{1},\ldots,X_{m}$ and $X_{m+1},\ldots X_{N}$
are freely independent and each generates a diffuse von Neumann algebra,
then for $T\in A\otimes A$ the derivation $\delta$ defined by $\delta(X_{j})=[X_{j},T]$
for $1\leq j\leq m$ and $\delta(X_{j})=0$ for $m+1\leq j\leq N$
is also in $\operatorname{Der}_{c}(A)$. Then one easily gets that
$\dim_{M\bar{\otimes}M^{o}}\overline{V}>1$ (indeed, $V$ contains
vectors of the form $([T,X_{1}],\ldots,[T,X_{m}],0,\ldots,0)$, $T\in L^{2}(A\otimes A)$,
and so its closure is strictly larger than the closure of the set
of all vectors $([T,X_{1}],\ldots,[T,X_{N}])$, $T\in L^{2}(A\otimes A)$).
Thus $\delta_{0}(X_{1},\ldots,X_{N})>1$ if $W^{*}(A)$ is $R^{\omega}$
embeddable.

If $X_{1},\ldots,X_{N}$ are such that their conjugate variables (see
\cite{dvv:entropy5}) are polynomials, then the difference quotient
derivations are in $\operatorname{Der}_{c}$ and thus $V=(A\otimes A)^{N}$,
and so $\delta_{0}=N$ (if $W^{*}(A)$ is $R^{\omega}$ embeddable).

In the case that $X_{1},\ldots,X_{N}$ are generators of the group
algebra $\mathbb{C}\Gamma$ of a discrete group $\Gamma$, $\delta^{*}(X_{1},\ldots,X_{N})=\beta_{1}^{(2)}(\Gamma)-\beta_{0}^{(2)}(\Gamma)+1$,
where $\beta_{j}^{(2)}$ are the $L^{2}$ Betti numbers of $\Gamma$
(see \cite{luck:book} for a definition). It is therefore natural
to ask whether the same holds true for $\delta_{0}$ instead of $\delta^{*}$
for some class of groups. If this is true, then knowing that $\beta_{1}^{(2)}(\Gamma)\neq0$
implies that $\delta_{0}>1$ and thus the group algebra has a variety
of properties that we explained above (see also \cite{peterson:l2rigid}).

It is clearly necessary for the equality $\delta_{0}=\beta_{1}^{(2)}-\beta_{0}^{(2)}+1$
that $\Gamma$ can be embedded into the ultrapower of the hyperfinite
II$_{1}$ factor (because otherwise $\delta_{0}$ would be $-\infty$).
In particular, one is tempted to conjecture that equality holds at
least in the case when $\Gamma$ is residually finite. 

Theorem \ref{thm:estAlg} implies a result like the one in \cite{brown-dykema-jung:fdimamalg}:

\begin{thm}
Assume that $\Gamma$ is embeddable into the unitary group of the
ultrapower of the hyperfinite II$_{1}$ factor. Then\[
\delta_{0}(\Gamma)\geq\dim_{L(\Gamma)}\overline{\{c:\Gamma\to\mathbb{C}\Gamma\textrm{ cocycle}\}}.\]
In particular, if $\Gamma$ belongs to the class of groups containing
all groups with $\beta_{1}^{(2)}=0$ and closed under amalgamated
free products over finite subgroups, passage to finite index subgroups
and finite extensions, then\[
\delta_{0}(\Gamma)=\beta_{1}^{(2)}(\Gamma)-\beta_{0}^{(2)}(\Gamma)+1.\]

\end{thm}
Let us now describe the main idea of the present paper. Our main result
states that if the von Neumann algebra $M=W^{*}(X_{1},\ldots,X_{n})$
can be embedded into the ultrapower of the hyperfinite II$_{1}$ factor,
then\begin{equation}
\delta_{0}(X_{1},\ldots,X_{n})\geq\dim_{M\bar{\otimes}M^{o}}V,\label{eq:mainineq}\end{equation}
where $V=\overline{\{(\partial(X_{1}),\dots,\partial(X_{n})):\partial\in\mathcal{C}\}}^{L^{2}}$
and $\mathcal{C}$ is some class of derivations from the algebra of
non-commutative polynomials $\mathbb{C}[X_{1},\ldots,X_{n}]$ to $L^{2}(M)\bar{\otimes}L^{2}(M^{o})$,
which will be made precise later.

The quantity $\delta_{0}(X_{1},\ldots,X_{n})$ is, very roughly, a
kind of Minkowski dimension ({}``relative'' to $R^{\omega}$) of
the set $\mathcal{V}$ of embeddings of $M$ into $R^{\omega}$, the
ultrapower of the hyperfinite II$_{1}$ factor (indeed, the set of
such embeddings can be identified with the set of images under the
embedding of the generators $X_{1},\ldots,X_{n}$, i.e., with the
set of microstates for $X_{1},\ldots,X_{n}$). On the other hand,
$\dim_{M\bar{\otimes}M^{o}}V$ is a linear dimension (relative to
$M\bar{\otimes}M^{o}$) of a certain vector space. If we could find
an interpretation for $V$ as a subspace of a {}``tangent space''
to $\mathcal{V}$, then the inequality \eqref{eq:mainineq} takes
the form of the inequality linking the Minkowski dimension of a manifold
with the linear dimension of its tangent space. One natural proof
of such an inequality would involve proving that a linear homomorphism
of the tangent space to a manifold at some point can be exponentiated
to a local diffeomorphism of a neighborhood of that point. 

Thus an essential step in proving a lower inequality on free entropy
dimension is to find an analog of such an exponential map. 

This leads to the idea, given a matrix $Q_{ij}\in(L^{2}(M)\bar{\otimes}L^{2}(M^{o}))^{n}$
of values of derivations (so that $Q_{ij}=\partial_{j}(X_{i})$ for
some $n$-tuple of derivation $\partial_{j}$ belonging to our class
$\mathcal{C}$), to try to associate to $Q$ a one-parameter deformation
$\alpha_{t}$ of a given embedding $\alpha=\alpha_{0}$ of $M$ into
$R^{\omega}$. It turns out that there are two (related) ways to do
this. 

The first approach comes from the idea that we (at least in principle)
know how to exponentiate derivations from an algebra to itself (the
result should be a one-parameter automorphism group of the algebra).
We thus try to extend $\partial=\partial_{1}\oplus\cdots\oplus\partial_{n}$
to a derivation of a larger algebra $\mathcal{A}=\mathbb{C}[X_{1},\ldots,X_{n},S_{1},\ldots,S_{n}]$,
where $S_{1},\ldots,S_{n}$ are free from $X_{1},\ldots,X_{n}$ and
form a free semicircular family. The key point is that the closure
in $L^{2}(\mathcal{A})$ of $MS_{1}M+\cdots+MS_{n}M$ is isomorphic
to $[L^{2}(M)\otimes L^{2}(M)]^{n}$. The inverse of this isomorphism
takes an $n$-tuple $a=(a_{1}\otimes b_{1},\ldots,a_{n}\otimes b_{n})$
to $\sum a_{j}S_{j}b_{j}$, which we denote by $a\#S$. We now define
a new derivation $\tilde{\partial}$ of $\mathcal{A}$ with values
in $L^{2}(\mathcal{A})$ by $\tilde{\partial}(X_{j})=\partial(X_{j})\#S$.
To be able to exponentiate $\tilde{\partial}$, we need to make sure
that it is anti-Hermitian as an unbounded operator on $L^{2}(\mathcal{A})$,
which naturally leads to the equation $\tilde{\partial}(S_{j})=-\partial^{*}(\zeta_{j})$,
where $\zeta_{j}=(0,\ldots,1\otimes1,\ldots,0)$ ($j$-th entry nonzero).
One can check that if $\zeta_{j}$ is in the domain of $\partial^{*}$
for all $j$, then $\tilde{\partial}$ is a closable operator which
has an anti-Hermitian extension, and so can be exponentiated to a
one-parameter group of automorphisms $\alpha_{t}$ of $L^{2}(\mathcal{A})$.
Unfortunately, unless we know more about the derivation $\partial$
(such as, for example, assuming that $\tilde{\partial}(\mathcal{A})\subset\mathcal{A}$),
we cannot prove that $\alpha_{t}$ takes $W^{*}(\mathcal{A})$ to
$W^{*}(\mathcal{A})$. However, if this is the case, then we do get
a one-parameter family of embeddings $\alpha_{t}|_{M}:M\to M*L(\mathbb{F}(n))\subset R^{\omega}$.
We explain this approach in more detail in the appendix Appendix (\S\ref{sec:owDer}).

The second approach was suggested to us by A. Guionnet, to whom we
are indebted for generously allowing us to publish it. The idea involves
considering the free stochastic differential equation\begin{equation}
dX_{j}(t)=\sum_{i}Q_{ij}(X_{1}(t),\ldots,X_{n}(t))\#dS_{i}-\frac{1}{2}\xi_{j}(X_{1}(t),\ldots,X_{n}(t)),\qquad X_{j}(0)=X_{j},\label{eq:freeSDE0}\end{equation}
where $\partial(X_{j})=(Q_{1j},\ldots,Q_{nj})\in(L^{2}(M)\bar{\otimes}L^{2}(M^{o}))^{n}$
and $\xi_{j}(X_{1},\ldots,X_{n})=\partial^{*}\partial(X_{j})$. One
difficulty in even phrasing the problem is that it is not quite clear
what is meant by $Q_{ij}$ and $\xi_{j}$ applied to their arguments
(in the classical case, this would mean a function applied to the
random variable $X(t)$). However, if this equation can be formulated
and has a stationary solution $X(t)$ (i.e., one for which the law
does not depend on $t$), then the map $\alpha_{t}:X_{j}\mapsto X_{j}(t^{2})$
determines a one-parameter family of embeddings of the von Neumann
algebra $M$ into some other von Neumann algebra $\mathcal{M}$ (generated
by all $X(t):t\geq0$). This can be carried out successfully if $Q$
and $\xi$ are sufficiently nice; this is this is the case, for example,
when $X_{1},\ldots,X_{n}$ are $q$-semicircular variables, in which
case $Q$ and $\xi$ can be taken to be analytic non-commutative power
series.

Let us assume now that $\partial$ takes $\mathcal{B}=\mathbb{C}[X_{1},\ldots,X_{n}]$
to $\mathcal{B}\otimes\mathcal{B}^{o}$ and also $\partial^{*}(1\otimes1)\in\mathcal{B}$\emph{
}(this is the case, for example, if $X_{1},\ldots,X_{n}$ have polynomial
conjugate variables \cite{dvv:entropy5}). Then both approaches work
to actually give one a stronger statement: one gets a one-parameter
family of embeddings $\alpha_{t}:M\to R^{\omega}$ so that $\Vert\alpha_{t}(X_{j})-(X_{j}+t\sum_{i}Q_{ij}S_{i})\Vert_{2}=O(t^{2})$.
Let us assume for the moment that $Q_{ij}=\delta_{ij}1\otimes1$,
so that our estimate reads\begin{equation}
\Vert\alpha_{t}(X_{j})-(X_{j}+tS_{j})\Vert_{2}=O(t^{2}).\label{eq:OWlikeEstimate0}\end{equation}
An estimate of this kind was used as a crucial step by Otto and Villani
in their work on the classical transportation cost inequality \cite[\S{}4 Lemma 2]{otto-villani:transportation};
a free version (for $n=1$) is the key ingredient in the proof of
free transportation cost inequality and free Wasserstein distance
given by Biane-Voiculescu \cite{bine-voiculescu:WassersteinDist}.
Indeed, since the law of $\alpha_{t}(X_{j})$ is the same as $X_{j}$,
one obtains after working out the error bounds an estimate on the
non-commutative Wasserstein distance between the laws $\mu_{X_{1},\ldots,X_{n}}$
and $\mu_{X_{1}+tS_{1},\ldots,X_{n}+tS_{n}}$:\[
d_{W}(\mu_{X_{1},\ldots,X_{n}},\mu_{X_{1}+tS_{1},\ldots,X_{n}+tS_{n}})\leq\frac{1}{2}\Phi(X_{1},\ldots,X_{n})^{1/2}t+O(t^{2}).\]
We now point out that this estimate is of direct relevance to a lower
estimate on $\delta_{0}$. Indeed, suppose that some $n$-tuple of
$k\times k$ matrices $x_{1},\ldots,x_{n}$ has as its law approximately
the law of $X_{1},\ldots,X_{n}$ (i.e., $(x_{1},\ldots,x_{n})\in\Gamma(X_{1},\ldots,X_{n};k,l,\varepsilon)$
in the notation of \cite{dvv:entropy2}). Then \eqref{eq:OWlikeEstimate0}
implies that by approximating $\alpha_{t}(X_{j})$ with polynomials
in $X_{1},\ldots,X_{n},S_{1},\ldots,S_{n}$, one can find another
$n$-tuple $x_{1}',\ldots,x_{n}'$ with almost the same law as $X_{1},\ldots,X_{n}$,
and so that $\Vert x_{j}'-(x_{j}+ts_{j})\Vert\leq Ct^{2}$ (here $s_{1},\ldots,s_{n}$
are some matrices whose law is approximately that of $S_{1},\ldots,S_{n}$,
and which are approximately free from $x_{1},\ldots,x_{n}$). But
this means that if one moves along a line starting at $x_{1},\ldots,x_{n}$
in the direction of $s_{1},\ldots,s_{n}$, then the distance to the
set $\Gamma(X_{1},\ldots,X_{n};k,l,\varepsilon)$ grows quadratically.
Thus this line is \emph{tangent} to the set $\Gamma(X_{1},\ldots,X_{n};k,l,\varepsilon)$.
From this one can derive estimates relating the packing numbers of
$\Gamma(X_{1},\ldots,X_{n};k,l,\varepsilon)$ and $\Gamma(X_{1}+tS_{1},\ldots,X_{n}+tS_{n};k,l,\varepsilon)$
which can be converted into a lower estimate on $\delta_{0}$.

In conclusion, it is worth pointing out that the main obstacle that
we face in trying to extend the estimate \eqref{eq:mainineq} to larger
classes of derivations is the question of existence of stationary
solutions of \eqref{eq:freeSDE0} for more general classes of functions
$Q$ and $\xi$ (and not, surprisingly enough, the {}``usual'' difficulties
in dealing with sets of microstates).

\begin{acknowledgement*}
The author is grateful to A. Guionnet for suggesting the idea of using
stationary solutions to free SDEs as an alternative form of {}``exponentiating''
derivations, and to (patiently) explaining to him about stochastic
differential equations. The author also thanks D.-V. Voiculescu for
a number of comments and suggestions.
\end{acknowledgement*}

\section{Existence of stationary solutions.}

\subsection{Free SDEs with analytic coefficients.}

The main result of this section states that a free stochastic differential
equation of the form\[
dX_{t}=\Xi\#dS_{t}-\frac{1}{2}\xi_{t}dt\]
where $X_{t}$ is an $N$-tuple of random variables has a stationary
solution, as long as the coefficients $\Xi$ and $\xi$ are analytic
(i.e., are non-commutative power series with sufficient radii of convergence).

\subsubsection{Estimates on certain operators appearing in free Ito calculus.\label{sub:PowerSeriesEst}}

Let $f$ be a non-commutative power series in $N$ variables. We will
denote by $c_{f}(n)$ the maximal modulus of a coefficient of a monomial
of degree $n$ in $f$. Thus if $f=\sum f_{i_{1}\dots i_{n}}X_{i_{1}}\cdots X_{i_{n}}$,
then $c_{f}(n)=\max_{i_{1}\dots i_{n}}|f_{i_{1}\dots i_{n}}|$. We'll
also write\[
\phi_{f}(z)=\sum c_{f}(n)z^{n}.\]
Then $\phi_{f}(z)$ is a formal power series in $z$. If $\rho$ is
the radius convergence of $\phi_{f}$, we'll say that $R=\rho/N$
is the multi-radius of convergence of $f$.

We'll also write\[
\Vert f\Vert_{\rho}=\sum_{n\geq0}c_{f}(n)N^{n}\rho^{n}\in[0,+\infty].\]
Note that $\Vert f\Vert_{\rho}=\sup_{|z|\leq N\rho}|\phi_{f}(z)|$
(since all of the coefficients in the power series $\phi_{f}(z)$
are real and positive).

We'll denote by $\mathcal{F}(R)$ the collection of all power series
$f$ for which the multi-radius of convergence is at least $R$. In
other words, we require $\Vert f\Vert_{\rho}<\infty$ for all $\rho<R$. 

Note that $\mathcal{F}_{R}$ is a complete topological vector space
if endowed with the topology $T_{i}\to T$ iff $\Vert T_{i}-T\Vert_{\rho}\to0$
for all $\rho<R$.

Let $\Psi$ be a non-commutative power series in $N$ variables having
the form\[
\sum f_{i_{1},\dots,i_{k};j_{1},\dots,j_{l}}Y_{i_{1}}\cdots Y_{i_{k}}\otimes Y_{j_{1}}\cdots Y_{j_{l}}.\]
We'll call $\Psi$ a formal non-commutative power series with values
in $\mathbb{C}\langle Y_{1},\dots,Y_{N}\rangle^{\otimes2}$. We'll
write $c_{\Psi}(m,n)$ the maximal modulus of a coefficient of a monomial
of the form $Y_{i_{1}}\cdots Y_{i_{m}}\otimes Y_{j_{1}}\cdots Y_{j_{n}}$
in $\Psi$. We let $\phi_{\Psi}(z,w)=\sum_{n,m}c_{\psi(m,n)}z^{m}w^{n}$.
We'll put\[
\Vert\Psi\Vert_{\rho}=\sup_{|z|,|w|\leq N\rho}|\phi_{\Psi}(z,w)|=\phi_{\Psi}(N\rho,N\rho)=\sum_{n\geq0}\left[\sum_{k+l=n}c_{\Psi}(k,l)\right]N^{n}\rho^{n}\in[0,+\infty].\]
We'll denote by $\mathcal{F}'(R)$ the collection of all non-commutative
power series for which $\Vert\Psi\Vert_{\rho}<\infty$ for all $\rho<R$. 

It will be convenient to use the following notation. Let $\phi(z_{1},\dots,z_{n})$,
$\psi(z_{1},\dots,z_{n})$ be two formal power series (in commuting
variables). We'll say that $\phi\prec\psi$ if all coefficients in
$\phi,\psi$ are real and positive, and for each $k_{1},\dots,k_{n}$,
the coefficient of $z_{1}^{k_{1}}\cdots z_{n}^{k_{n}}$ in $\phi$
is less than or equal to the corresponding coefficient in $\psi$. 

If \emph{$\mathcal{M}$ }is a unital Banach algebra, $Y_{1},\dots,Y_{N}\in\mathcal{M}$
and $\Vert Y_{j}\Vert<\rho$ for all $j$, then $\Vert g(Y_{1},\dots,Y_{n})\Vert\leq\Vert g\Vert_{\rho}$
whenever $g$ is in any one of the spaces $\mathcal{F}(R)$, or $\mathcal{F}'(R)$
(here the norm $\Vert g(Y_{1},\dots,Y_{n})\Vert$ denotes the norm
on $\mathcal{M}$ or on the projective tensor product $\mathcal{M}^{\otimes2}$,
as appropriate).

We now collect some facts about power series:

\begin{itemize}
\item Let $f,g\in\mathcal{F}(R)$. Then $\phi_{fg}\prec\phi_{f}\phi_{g}$.
In particular, $fg\in\mathcal{F}(R)$ and $\Vert fg\Vert_{\rho}\le\Vert f\Vert_{\rho}\Vert g\Vert_{\rho}$.
\item Let $f=\sum f_{i_{1}\dots i_{n}}X_{i_{1}}\cdots X_{i_{n}}\in\mathcal{F}(R)$
and denote by $\mathcal{D}_{ij}f$ the formal power series\[
\mathcal{D}_{ij}f=\sum_{i_{1}\dots i_{n}}\sum_{k<l}\delta_{i_{k}=i}\delta_{i_{l}=j}f_{i_{1}\dots i_{n}}X_{i_{k+1}}\cdots X_{i_{l-1}}\otimes X_{i_{l+1}}\cdots X_{i_{n}}X_{i_{1}}\cdots X_{i_{k-1}}.\]
Since a monomial $X_{i_{1}}\cdots X_{i_{k}}\otimes X_{j_{1}}\cdots X_{j_{r}}$
could arise in the expression for $\mathcal{D}_{ij}f$ in at most
$r+1$ ways, $c_{\mathcal{D}_{j}f}(a,b)\leq(b+1)c_{f}(a+b+2)$. Denote
by $\hat{\phi}_{f}$ the power series $\hat{\phi}_{f}(z,w)=\sum_{n,m}(n+1)c_{f}(n+m+2)z^{m}w^{n}$.
Then $\phi_{\mathcal{D}_{j}f}\prec\hat{\phi}_{f}$. Since $\hat{\phi}_{f}(z,z)\prec\phi_{f}''(z)$,
we conclude that\[
\Vert\mathcal{D}_{ij}f\Vert_{\rho}\leq\sup_{|z|\leq N\rho}|\phi_{f}''(z)|\]
and in particular $\mathcal{D}_{ij}f\in\mathcal{F}'(R)$. 
\item Let $\Theta=\sum\Theta_{i_{1}\dots i_{n};j_{1}\dots j_{m};k_{1}\dots k_{p}}X_{i_{1}}\cdots X_{i_{n}}\otimes X_{j_{1}}\cdots X_{j_{m}}\in\mathcal{F}'(R)$,
and let $\Psi=\sum\Psi_{i_{1}\dots i_{n};j_{1}\dots j_{m}}X_{i_{1}}\dots X_{i_{n}}\otimes X_{j_{1}}\dots X_{j_{m}}\in\mathcal{F}'$.
Consider\[
\Psi\#_{in}\Theta=\sum\Psi_{t_{1}\dots t_{a},s_{1},\dots,s_{b}}\Theta_{i_{1}\dots i_{n};j_{1}\dots j_{m}}X_{i_{1}}\cdots X_{i_{n}}\ X_{t_{1}}\cdots X_{t_{a}}\otimes X_{s_{1}}\cdots X_{s_{b}}\ X_{j_{1}}\cdots X_{j_{m}}.\]
(In the simple case that $\Psi=A\otimes B$ and $\Theta=P\otimes Q$,
where $A,B,P,Q$ are monomials, $\Psi\#_{in}\Theta=PA\otimes BQ$,
i.e., $\#_{in}$ is the {}``inside'' multiplication on \emph{$\mathcal{F}'(R)$}).
Then\[
c_{\Psi\#_{in}\Theta}(n,m)\leq\sum_{k+l=n}\sum_{r+s=m}c_{\Psi}(k,r)c_{\Theta}(l,s),\]
so the coefficient of $z^{n}w^{m}$ in $\phi_{\Psi\#_{in}\Theta}(z,w)$
is dominated by the coefficient of $z^{n}w^{m}$ in $\phi_{\Psi}(z,w)\phi_{\Theta}(z,w)$.
Consequently, $\phi_{\Psi\#_{in}\Theta}\prec\phi_{\Psi}\phi_{\Theta}$
and\[
\Vert\Psi\#_{in}\Theta\Vert_{\rho}\leq\Vert\Psi\Vert_{\rho}\Vert\Theta\Vert_{\rho}.\]
In particular, $\Psi\#_{in}\Theta\in\mathcal{F}'(R)$. Similar estimates
and conclusion of course hold for the {}``outside'' multiplication
$\Psi\#_{out}\Theta$, defined by\[
\Psi\#_{out}\Theta=\sum\Psi_{s_{1},\dots,s_{b};t_{1}\dots t_{a}}\Theta_{i_{1}\dots i_{n};j_{1}\dots j_{m}}X_{t_{1}}\cdots X_{t_{a}}\ X_{i_{1}}\cdots X_{i_{n}}\otimes X_{j_{1}}\cdots X_{j_{m}}\ X_{s_{1}}\cdots X_{s_{b}}.\]
In that case we get $\phi_{\Psi\#_{out}\Theta}(z,w)\prec\phi_{\Psi}(w,z)\phi_{\Theta}(z,w)$
and $\Vert\Psi\#_{out}\Theta\Vert_{\rho}\leq\Vert\Psi\Vert_{\rho}\Vert\Theta\Vert_{\rho}$.
\item Let $\tau$ be a linear functional on the algebra of non-commutative
polynomials in $n$ variables, so that $|\tau(X_{i_{1}}\cdots X_{i_{n}})|\le R_{0}^{n}$
for all $n$. Given $\Theta=\sum\Theta_{i_{1}\dots i_{n};j_{1}\dots j_{m}}X_{i_{1}}\cdots X_{i_{n}}\otimes X_{j_{1}}\cdots X_{j_{m}}\in\mathcal{F}'(R)$,
assume that $R_{0}<R$ and consider the formal sum\[
(1\otimes\tau)(\Theta)=\sum_{n,i_{1},\dots,i_{n}}\left[\sum_{m,j_{1},\dots,j_{m}}\Theta_{i_{1}\dots i_{n};j_{1}\dots j_{m}}\tau(X_{j_{1}}\cdots X_{j_{m}})\right]X_{i_{1}}\cdots X_{i_{n}}.\]
More precisely, we consider the formal power series in which the coefficient
of $X_{i_{1}}\cdots X_{i_{n}}$ is given by the sum\[
\sum_{m,j_{1},\dots,j_{m}}\Theta_{i_{1}\dots i_{n};j_{1}\dots j_{m}}\tau(X_{j_{1}}\cdots X_{j_{m}}).\]
But since $|\tau(X_{j_{1}}\cdots X_{j_{m}})|\leq R_{0}^{m}$, this
sum is bounded by the coefficient of $z^{n}$ in the power series
expansion of $\phi(z,NR_{0})$ (as a function of $z$), and is convergent.
Thus $\phi_{(1\otimes\tau)(\Theta)}(z)\prec\phi_{\Theta}(z,NR_{0})$
and we readily see that $(1\otimes\tau)(\Theta)$ is well-defined,
belongs to $\mathcal{F}(R)$ and moreover\[
\Vert(1\otimes\tau)(\Theta)\Vert_{\rho}\leq\Vert\Theta\Vert_{\rho}\]
whenever $\rho>R_{0}$. 
\end{itemize}
We now combine these estimates:

\begin{lem}
\label{lemma:estimatesQ}Let $\tau$ as above be a linear functional
on the space of non-commutative polynomials in $N$ variables satisfying
$\tau(X_{i_{1}}\cdots X_{i_{n}})\leq R_{0}^{n}$. Let $R>R_{0}$ and
assume that $\xi_{j}\in\mathcal{F}(R)$, $j=1,\dots,N$, $\Psi=(\Psi_{ij})\in M_{N\times N}\mathcal{F}'(R)$.
For $f\in\mathcal{F}(R)$ let\[
\mathcal{L}^{(\tau)}(f)=(1\otimes\tau)(\sum_{ijk}\Psi_{jk}\#_{in}(\Psi_{ki}\#_{out}(\mathcal{D}_{ij}f)))-\sum_{j}\frac{1}{2}\xi_{j}f.\]
Then $\mathcal{L}_{j}^{(\tau)}(f)\in\mathcal{F}(R)$ and moreover
for any $R_{0}<\rho<R$,\begin{eqnarray*}
\Vert\mathcal{L}^{(\tau)}(f)\Vert_{\rho} & \leq & \sum_{ijk}\Vert\Psi_{jk}\Vert_{\rho}\Vert\Psi_{ki}\Vert_{\rho}\cdot\sup_{|z|\leq N\rho}|\phi''_{f}(z)|+\frac{1}{2}\sum_{j}\Vert\xi_{j}\Vert_{\rho}\Vert f\Vert_{\rho}\\
\phi_{\mathcal{L}^{(\tau)}(f)}(z) & \prec & \sum_{ijk}\phi_{\Psi_{jk}}(z,NR_{0})\phi_{\Psi_{ki}}(NR_{0},z)\hat{\phi}_{f}(z,NR_{0})+\frac{1}{2}\sum_{j}\phi_{\xi_{j}}(z)\phi_{f}(z).\end{eqnarray*}
where $\hat{\phi}_{f}(z,w)=\sum_{n,m}(n+1)c_{f}(n+m+2)z^{m}w^{n}$. 
\end{lem}
For $\phi$ a power series in $z,w_{1},\dots,w_{k}$ with multi-radius
of convergence bigger than $\rho$ and all coefficients of monomials
non-negative, let $\phi_{w_{1},\dots,w_{k}}(z)=\phi(z,w_{1},\dots,w_{k})$.
Set\begin{eqnarray*}
Q\phi(z,w_{1},\dots w_{k+1}) & = & \widehat{\phi_{w_{1},\dots,w_{k}}}(z,w_{k+1})\\
D\phi(z,w_{1},\dots,w_{k}) & = & \partial_{z}^{2}\phi(z,w_{1},\dots,w_{k}).\end{eqnarray*}

We note that $\hat{\phi}(z,z)\prec\phi''(z)$, and that both $Q$
and $D$ are monotone for the ordering $\prec$. It follows that if
$\kappa_{j}$, $\lambda_{j}$ are some power series with radius of
convergence bigger than $\rho$ and positive coefficients, then for
any $a_{1},b_{1},\dots,a_{k},b_{k}\geq0$ and any $R<\rho$,\begin{multline*}
\left[Q^{a_{1}}\kappa_{1}(z)D^{b_{1}}\lambda_{1}(z)Q^{a_{2}}\kappa_{2}(z)D^{b_{2}}\lambda_{2}(z)\cdots D^{b_{k}}\lambda_{k}\right]\Big|_{z=w_{1}=\cdots=w_{\sum b_{k}}=R}\\
\leq\left[D^{a_{1}}\kappa_{1}(z)D^{b_{1}}\lambda_{1}(z)D^{a_{2}}\kappa_{2}(z)D^{b_{2}}\lambda_{2}(z)\cdots D^{b_{k}}\lambda_{k}\right]\Big|_{z=w_{1}=\cdots=w_{\sum b_{k}}=R}\end{multline*}
Define now\[
\hat{\mathcal{L}}\phi(z)=\sum_{ijk}\phi_{\psi_{jk}}(z,NR_{0})\phi_{\Psi_{ki}}(NR_{0},z)\phi''(z)+\frac{1}{2}\sum_{j}\phi_{\xi_{j}}(z)\phi(z).\]
Then we have obtained the following inequality:\[
\phi_{\mathcal{L}^{n}f}(NR_{0})\leq\hat{\mathcal{L}}^{n}\phi_{f}(NR_{0}).\]
We record this observation:

\begin{lem}
\label{lemma:estimateL}Let $\hat{\mathcal{L}}\phi(z)=\sum_{ijk}\phi_{\psi_{jk}}(z,NR_{0})\phi_{\Psi_{ki}}(NR_{0},z)\phi''(z)+\frac{1}{2}\sum_{j}\phi_{\xi_{j}}(z)\phi(z)$
and let $\tau$ be a trace so that for any monomial $P$, $|\tau(P)|<R_{0}^{n}$,
$n=\deg P$. Then\[
|\tau(\mathcal{L}^{n}f)|\leq\hat{\mathcal{L}}^{n}\phi_{f}(NR_{0}).\]

\end{lem}

\subsubsection{Analyticity of $\partial^{*}\partial(X_{j})$.}

Let us now assume that $\Xi=(\Xi_{1},\dots,\Xi_{N})\in\mathcal{F}'(R)$.
Let $(X_{1},\dots,X_{N})$ be an $N$-tuple of self-adjoint operators
in a tracial von Neumann algebra $(M,\tau)$, and assume that $\Vert X_{j}\Vert<R$
for all $j$. Let $\partial:L^{2}(M)\to L^{2}(M)\bar{\otimes}L^{2}(M)$
be the derivation densely defined on polynomials in $X_{1},\dots,X_{N}$
by $\partial(X_{j})=\Xi_{j}(X_{1},\dots,X_{N})$. We'll assume that
$1\otimes1$ belongs to the domain of $\partial^{*}$ and that there
exists some $\zeta\in\mathcal{F}(R)$ so that $\partial^{*}(1\otimes1)=\zeta(X_{1},\dots,X_{N})$.

\begin{lem}
\label{lem:1otimes1Enough}With the above assumptions, there exist
$\xi_{j}\in\mathcal{F}(R)$, $j=1,\dots,N$, so that $\xi_{j}(X_{1},\dots,X_{N})=\partial^{*}\partial(X_{j})$.
\end{lem}
\begin{proof}
It follows from \cite{dvv:entropy5,shlyakht:amalg} that under these
assumptions, $\partial$ is closable. Moreover, for any $a,b$ polynomials
in $X_{1},\dots,X_{N}$, $a\otimes b$ belongs to the domain of $\partial^{*}$
and \[
\partial^{*}(a\otimes b)=a\zeta b+(1\otimes\tau)[\partial(a)]b+a(\tau\otimes1)[\partial(b)],\]
where $\zeta=\zeta(X_{1},\dots,X_{N})=\partial^{*}(1\otimes1)$.

Consider now formal power series in $N$ variables having the form\[
\Theta=\sum\Theta_{i_{1},\dots,i_{k};j_{1},\dots,j_{l};t_{1},\dots,t_{r}}Y_{i_{1}}\cdots Y_{i_{k}}\otimes Y_{j_{1}}\cdots Y_{j_{l}}\otimes Y_{t_{1}}\cdots Y_{t_{r}}.\]
We'll write $\phi_{\Theta}(z,w,v)$ for the power series whose coefficient
of $z^{m}w^{n}v^{k}$ is equal to the maximum\[
\max\{|\Theta_{i_{1},\dots,i_{m};j_{1},\dots,j_{n};t_{1},\dots,t_{k}}|:i_{1},\dots,i_{m},j_{1},\dots,j_{n},t_{1},\dots,t_{n}\in\{1,\dots,N\}\}.\]
We'll denote by $\mathcal{F}''(R)$ the collection of all such power
series for which $\phi_{\Theta}$ has a multi-radius of convergence
at least $NR$. 

Let $\mathcal{D}_{1}^{(s)}:\mathcal{F}'(R)\to\mathcal{F}''(R)^{}$
be given by\begin{multline*}
\mathcal{D}_{1}^{(s)}\sum f_{i_{1},\dots,i_{k};j_{1},\dots,j_{l}}Y_{i_{1}}\cdots Y_{i_{k}}\otimes Y_{j_{1}}\cdots Y_{j_{l}}=\\
\sum f_{i_{1},\dots,i_{k};j_{1},\dots,j_{l}}\sum_{p}\delta_{i_{p}=s}Y_{i_{1}}\cdots Y_{i_{p-1}}\otimes Y_{i_{p+1}}\cdots Y_{i_{k}}\otimes Y_{j_{1}}\cdots Y_{j_{l}}.\end{multline*}
Then clearly $\phi_{\mathcal{D}_{1}^{(s)}(\Psi)}(z,z,w)\prec\partial_{z}\phi_{\Psi}(z,w)$,
so that $\mathcal{D}_{1}^{(s)}\Psi$ indeed lies in $\mathcal{F}''(R)$
if $\Psi\in\mathcal{F}'(R)$. 

Similarly, if we define for $\Psi\in\mathcal{F}'(R)$, $\Theta\in\mathcal{F}''(R)$\[
\Psi\#_{in}^{(1)}\Theta=\sum\Psi_{t_{1}\dots t_{a},s_{1},\dots,s_{b}}\Theta_{i_{1}\dots i_{n};j_{1}\dots j_{m};k_{1}\dots k_{p}}Y_{i_{1}}\cdots Y_{i_{n}}\ Y_{t_{1}}\cdots Y_{t_{a}}\otimes Y_{s_{1}}\cdots Y_{s_{b}}\ Y_{j_{1}}\cdots Y_{j_{m}}\otimes Y_{k_{1}}\cdots Y_{k_{p}},\]
then $\phi_{\Psi\#_{in}^{(1)}\Theta}(z,v,w)\prec\phi_{\Psi}(z,v)\phi_{\Theta}(z,v,w)$
and in particular $\Psi\#_{in}^{(1)}\Theta\in\mathcal{F}''(R)$. (Note
that $\#_{in}^{(1)}$ corresponds to {}``multiplying around'' the
first tensor sign in $\Theta$).

Finally, if $\tau$ is any linear functional so that $\tau(P)<R_{0}^{\deg P}$
for any monomial $P$ and we put\[
M_{2}(\Psi)=\sum\Psi_{i_{1},\dots,i_{n};j_{1},\dots,j_{m};k_{1},\dots k_{p}}Y_{i_{1}}\cdots Y_{i_{n}}\tau(Y_{j_{1}}\cdots Y_{j_{m}}Y_{k_{1}}\cdots Y_{k_{p}})\]
then $\phi_{M_{2}(\Psi)}(z)\leq\phi_{\Psi}(z,NR_{0},NR_{0})$ and
in particular $M_{2}(\Psi)\in\mathcal{F}(R)$ once $\Psi\in\mathcal{F}''(R)$
and $R_{0}<R$. In the foregoing, we'll use the trace $\tau$ of $M$
as our functional.

So if we put\[
T_{1}\Theta=M_{2}(\sum_{s}\Xi_{s}\#_{in}^{(1)}\mathcal{D}_{1}^{(s)}),\]
then $T_{1}$ maps $\mathcal{F}'(R)$ into $\mathcal{F}(R)$.

Note that in the case that $\Theta=A\otimes B$, where $A,B$ are
monomials, $T_{1}\Theta=(1\otimes\tau)(\partial(A)B)$. 

One can similarly define $T_{2}:\mathcal{F}'(R)\to\mathcal{F}(R)$;
it will have the property that $T_{2}\Theta=(1\otimes\tau)(\partial(A)B)$. 

Lastly, let $\zeta\in\mathcal{F}(R)$ and let $m:\mathcal{F}'(R)\to\mathcal{F}(R)$
be given by\[
m(\Theta)=\sum\Theta_{i_{1},\dots,i_{n};j_{1},\dots,j_{m}}\zeta_{p_{1},\dots,p_{r}}Y_{i_{1}}\cdots Y_{i_{n}}Y_{p_{1},\dots p_{r}}Y_{j_{1}}\cdots Y_{j_{m}}.\]
Once again, $\phi_{m(\Theta)}(z)\prec\phi_{\Theta}(z,z)\phi_{\zeta}(z)$.

Let now $Q(\Xi)=T_{1}(\Xi)+T_{2}(\Xi)+m(\Xi)$. We claim that $\xi=(Q(\Xi))(X_{1},\dots,X_{N})=\partial^{*}(\Xi(X_{1},\dots,X_{N}))$. 

Note that if $\Xi_{n}$ is a partial sum of $\Xi$ (say obtained as
the sum of all monomials in $\Xi$ having degree at most $n$), then
$Q(\Xi_{n})(X_{1},\dots,X_{N})=\partial^{*}(\Xi_{n}(X_{1},\dots,X_{N}))$.
Moreover, $ $as $n\to\infty$, $\Xi_{n}(X_{1},\dots,X_{N})\to\Xi(X_{1},\dots,X_{N})$
in $L^{2}$ and also $Q(\Xi_{n})(X_{1},\dots,X_{N})\to Q(\Xi)(X_{1},\dots,X_{N})$
in $L^{2}$ (this can be seen by observing first that the coefficients
of $Q_{n}(\Xi)$ converge to the coefficients of $Q(\Xi)$ and then
approximating $Q(\Xi)$ and $Q(\Xi_{n})$ by their partial sums). 

Since $\partial^{*}$ is closed, the claimed equality follows. 
\end{proof}

\subsubsection{Existence of solutions.}

Recall that a process $X_{1}^{(t)},\dots,X_{N}^{(t)}\in(M,\tau)$
is called stationary if its law does not depend on $t$; in other
words, for any polynomial $f$ in $N$ non-commuting variables, $\tau(f(X_{1}^{(t)},\dots,X_{N}^{(t)}))$
is constant.

\begin{lem}
\label{lem:SDEstationary} Let $X_{1}^{(0)},\dots,X_{N}^{(0)}$ be
an $n$-tuple of non-commutative random variables, $R_{0}>\max_{j}\Vert X_{j}^{(0)}\Vert$
and $R>R_{0}$. Let $\xi_{j}\in\mathcal{F}(R)$, $\Psi=(\Psi_{ij})\in M_{N\times N}(\mathcal{F}'(R))$,
so that $\Psi_{ij}(Z_{1},\dots,Z_{N})^{*}=\Psi_{ji}(Z_{1},\dots,Z_{N})$
for any self-adjoint $Z_{1},\dots,Z_{N}$.

Consider the free stochastic differential equation\begin{equation}
dX_{i}(t)=\Psi(X_{1}(t),\ldots,X_{N}(t))\#(dS_{t}^{(1)},\dots,dS_{t}^{(N)})-\frac{1}{2}\xi_{i}(X_{1}(t),\ldots,X_{N}(t))dt\label{eq:SDE}\end{equation}
with the initial condition $X_{j}(0)=X_{j}^{(0)}$, $j=1,\dots,n$.
Here $dS_{t}^{(1)},\ldots,dS_{t}^{(N)}$ is free Brownian motion,
and for $Q_{kl}=\sum a_{i}^{kl}\otimes b_{i}^{kl}\in M\hat{\otimes}M$,
and $Q=(Q_{kl})\in M_{N\times N}(M\hat{\otimes}M)$, we write $Q\#(W_{1},\dots,W_{N})=(\sum_{ki}a_{i}^{1k}W_{k}b_{i}^{1k},\dots,\sum a_{i}^{Nk}Wb_{i}^{Nk})$. 

Let $A=W^{*}(X_{1}^{(0)},\dots,X_{N}^{(0)})$ and let $\partial_{j}:L^{2}(A)\to L^{2}(A\bar{\otimes}A)$
be derivations densely defined on polynomials in $X_{1}^{(0)},\dots,X_{N}^{(0)}$
and determined by\[
\partial_{j}(X_{i})=Q_{ji}(X_{1}^{(0)},\dots,X_{N}^{(0)}).\]
Assume that for all $j$, $\partial_{i}X_{j}\in\operatorname{domain}\partial_{i}^{*}$
and that \[
\xi_{j}(X_{1}^{(0)},\dots,X_{N}^{(0)})=\sum_{i}\partial_{i}^{*}\partial_{i}(X_{j}^{(0)}).\]

Then there exists a $t_{0}>0$ and a stationary solution $X_{j}(t)$,
$0\leq t<t_{0}$. This stationary solution satisfies $X_{j}(t)\in W^{*}(X_{1},\ldots,X_{N},\{S_{j}(s):0\leq s\leq t\}_{j=1}^{N})$.
\end{lem}
We note that in view of Lemma \ref{lem:1otimes1Enough}, we may instead
assume that $1\otimes1\in\operatorname{domain}\partial_{j}^{*}$ and
$\partial_{j}^{*}(1\otimes1)=\zeta_{j}(X_{1}^{(0)},\dots,X_{N}^{(0)})$
for some $\zeta_{1},\dots,\zeta_{N}\in\mathcal{F}(R)$, since this
assumption guarantees the existence of $\xi_{j}\in\mathcal{F}(R)$
satisfying the hypothesis of Lemma \ref{lem:SDEstationary}. 

\begin{proof}
We note that, because $\Psi$ and $\xi$ are analytic, they are (locally) Lipschitz
in their arguments.

Thus it follows from the standard Picard argument (cf. \cite{biane-speicher:stochcalc})
that a solution (with given initial conditions) exists, at least for
all values of $t$ lying in some small interval $[0,t_{0})$, $t_{0}>0$.
Choose now $t_{0}$ so that $\Vert X_{j}(t)\Vert_{\infty}\leq R_{0}<R$
for all $0\leq t<t_{0}$ (this is possible, since the solution to
the SDE is locally norm-bounded).

Next, we note that if we adopt the notations of Lemma \ref{lemma:estimatesQ}
and define for $f\in\mathcal{F}(R)$\[
\mathcal{L}^{(\tau)}(f)=\sum_{ijk}(1\otimes\tau)(\Psi_{jk}\#_{in}(\Psi_{ki}\#_{out}(\mathcal{D}_{ij}f)))-\frac{1}{2}\sum_{j}\xi_{j}f,\]
then we have that $\mathcal{L}^{(\tau_{t})}f\in\mathcal{F}(R)$ (here
$\tau_{t}$ refers to the trace on $\mathbb{C}\langle X_{1}(t),\dots,X_{n}(t)\rangle$
obtained by restricting the trace from the von Neumann algebra containing
the process $X_{t}$ for small values of $t$, i.e., $\tau_{t}(P)=\tau(P(X_{1}(t),\dots,X_{n}(t)))$).
Ito calculus shows that for any $f\in\mathcal{F}(R)$,\[
\frac{d}{dt}\tau(f(X_{1}(t),\ldots,X_{N}(t)))\Big|_{t=s}=\tau_{s}\left((\mathcal{L}^{(\tau_{s})}f)(X_{1}(s),\ldots,X_{N}(s))\right).\]
In particular, replacing $f$ with $\mathcal{L}^{(\tau_{t})}f$ and
iterating gives us the equality\[
\left(\frac{d^{n}}{dt^{n}}\right)\tau(f(X_{1}(t),\ldots,X_{N}(t)))\Big|_{t=s}=\tau_{s}\left(((\mathcal{L}^{(\tau_{s})})^{n}f)(X_{1}(s),\ldots,X_{N}(s))\right).\]
Since $\xi_{j}(X_{1}(0),\ldots,X_{n}(0))=\sum_{i}\partial_{i}^{*}\partial_{i}(X_{j}(0))$,\[
\mathcal{L}^{(\tau_{0})}(f(X_{1}(0),\ldots,X_{N}(0)))=0\]
for any $f\in\mathcal{F}(R)$. Applying this to $f$ replaced with
$\mathcal{L}^{(\tau_{0})}f$ and iterating allows us to conclude that\[
\frac{d^{n}}{dt^{n}}\tau(f(X_{1}(t),\ldots,X_{N}(t)))\Big|_{t=0}=0,\qquad n\geq1.\]

Let\[
C_{n}(f,t)=\sup_{0\leq s\leq t}\Vert(\mathcal{L}^{(\tau_{s})})^{n}f(X_{1}(s),\ldots,X_{N}(s))\Vert.\]

By Lemma \ref{lemma:estimateL}, we have that\[
C_{n}(f,t)\leq|\hat{\mathcal{L}}^{n}\phi(R_{0})|,\]
where $\phi=\phi_{f}$ and\[
\hat{\mathcal{L}}\phi(z)=\sum_{ijk}\phi_{\Psi_{ik}}(z,NR_{0})\phi_{\Psi_{jk}}(NR_{0},z)\phi''(z)+\frac{1}{2}\sum_{j}\phi_{\xi_{j}}(z)\phi(z).\]
Thus if we set $C_{n}=|\hat{\mathcal{L}}^{n}\phi(R_{0})|$, then (because
all derivatives at zero of the function $\tau(f(X_{1}(t),\dots,X_{N}(t)))$
vanish),\begin{multline*}
|\tau(f(X_{1}(t),\dots,X_{N}(t)))-\tau(f(X_{1}(0),\dots,f(X_{N}(0)))|\\
=\left|\int_{0}^{t}\cdots\int_{0}^{t}\frac{d^{n}}{dr^{n}}\tau(f(X_{1}(r),\dots,f(X_{N}(r)))\Big|_{r=s}\ (ds)^{n}\right|\\
\leq\int_{0}^{t}\cdots\int_{0}^{t}C_{n}\ (ds)^{n}\leq C_{n}\frac{t^{n}}{n!}.\end{multline*}

We now note that $\hat{\mathcal{L}}\phi=\alpha_{1}(z)\phi''(z)+\alpha_{2}(z)\phi(z)$,
for some $\alpha_{1},\alpha_{2}$ analytic on $|z|<R$. Let $Q_{1}\phi=\alpha_{1}(z)\phi''$,
$Q_{2}\phi=\alpha_{2}\phi$, so that $\hat{\mathcal{L}}=Q_{1}+Q_{2}$.
Then if $\phi$ is a power series with all coefficients non-negative,
then so are $Q_{i}\phi$, $i=1,2$. Moreover, for any $\theta$ with
non-negative coefficients, $\theta Q_{1}\phi\prec Q_{1}(\theta\phi)$;
in particular, $Q_{2}Q_{1}\phi\prec Q_{1}Q_{2}\phi$. Furthermore,
if $\phi\prec\theta$ then $Q_{i}\phi\prec Q_{i}\theta$, $i=1,2$.
It follows that if $i_{1},\dots,i_{n}\in\{1,2\}$ are arbitrary and
exactly $k$ of $i_{1},\dots,i_{n}$ equal $1$, then\[
Q_{i_{1}}\cdots Q_{i_{n}}\phi\prec Q_{1}^{k}Q_{2}^{n-k}\phi=\frac{d^{k}}{dz^{k}}(\alpha_{i_{1}}\cdots\alpha_{i_{n}}\phi).\]
Thus for any $\rho\in(R_{0},R)$ and any $\phi$ with non-negative
coefficients, and assuming that $\rho-R_{0}<1$,\begin{eqnarray*}
Q_{i_{1}}\cdots Q_{i_{n}}\phi(R_{0}) & = & \frac{k!}{2\pi i}\int_{|w|=\rho}\frac{\alpha_{i_{1}}(w)\cdots\alpha_{i_{n}}(w)\phi(w)}{(w-R_{0})^{k+1}}dw\\
 & \leq & \frac{K^{n}Cn!}{(\rho-R_{0})^{n}},\end{eqnarray*}
where $K=\sup\{|\alpha_{i}(z)|:|z|=\rho,i=1,2\}$, $C=(2\pi(\rho-R_{0}))^{-1}\sup\{|\phi(z)|:|z|=\rho\}$. 

Since $\hat{\mathcal{L}}^{n}=\sum_{i_{1},\dots,i_{n}\in\{1,2\}}Q_{i_{1}}\cdots Q_{i_{n}}$,
we conclude that\[
C_{n}=\hat{\mathcal{L}}^{n}\phi(R_{0})\leq\sum_{i_{1},\dots,i_{n}\in\{1,2\}}Q_{i_{1}}\cdots Q_{i_{n}}\phi(R_{0})\leq Cn!\left[\frac{2K}{\rho-R_{0}}\right]^{n}.\]

Thus\[
|\tau(f(X_{1}(t),\dots,X_{N}(t)))-\tau(f(X_{1}(0),\dots,f(X_{N}(0)))|\leq C\left[\frac{2Kt}{\rho-R_{0}}\right]^{n}.\]
Thus we may choose $t_{0}$ small enough so that for any $t<t_{0}$,
$\left[2Kt(\rho-R_{0})^{-1}\right]^{n}\to0$ and so the solution is
indeed stationary on this interval.
\end{proof}
We note that once the equation \eqref{eq:SDE} has a stationary solution
on a small interval $[0,t_{0})$, then it of course has a stationary
solution for all time (since the same lemma applied to $X_{t_{0}/2}$
guarantees existence of the solution for up to $3t_{0}/2$ and so
on). However, we will not need this here.

\section{Otto-Villani type estimates.}

The main result of this section is an estimate on the non-commutative
Biane-Voiculescu-Wasserstein distance between the law of an $N$-tuple
of variables $X=X_{1},\ldots,X_{N}$ and the law of the $N$-tuple
$X+\sqrt{t}Q\#S$, where $S=(S_{1},\ldots,S_{N})$ is a free semicircular
family, $Q\in M_{N\times N}(L^{2}(W^{*}(X_{1},\ldots,X_{N})^{\otimes2}))$
is a matrix, and for $Q_{ij}=\sum_{k}A_{ij}^{(k)}\otimes B_{ij}^{(k)}$,
we denote by $Q\#S$ the $N$-tuple $(Y_{1},\ldots,Y_{N})$ with\[
Y_{i}=\sum_{j}\sum_{k}A_{ij}^{(k)}S_{j}B_{ij}^{(k)}.\]
Note that the sum defining $Y_{i}$ is operator-norm convergent; in
fact, the operator norm of $Y_{i}$ is the same as the $L^{2}$ norm
of the element\[
\sum_{j}\sum_{k}A_{ij}^{(k)}\otimes B_{ij}^{(k)}.\]

The estimate on Wasserstein distance (Proposition \ref{prop:OWprocess})
is obtained under the assumptions that a certain derivation, defined
by $\partial(X_{i})=(Q_{i1},\ldots,Q_{iN})\in(L^{2}(W^{*}(X_{1},\ldots,X_{N})^{\otimes2})^{N}$
is closable and satisfies certain further analyticity conditions (see
below for more precise statements). Under such assumptions, the estimate
states that\[
d_{W}(X,X+\sqrt{t}Q\#S)\leq Ct.\]
The main use of this estimate will be to give a lower bound for the
microstates free entropy dimension of $X_{1},\ldots,X_{N}$ (see Section
\ref{sec:deltaEst}).

\subsection{An Otto-Villani type estimate on Wasserstein distance via free SDEs.}

\begin{prop}
\label{prop:OWprocess}Let $\Xi\in M_{N\times N}(\mathcal{F}'(R))$,
$M=W^{*}(X_{1},\dots,X_{N}^{})$ and let $\partial_{j}:L^{2}(M)\to L^{2}(M\bar{\otimes}M)$
be derivations densely defined on polynomials in $X_{1},\dots,X_{N}$
and determined by\[
\partial_{j}(X_{i})=\Xi_{ji}(X_{1},\dots,X_{N}).\]
Assume that for all $j$, $1\otimes1\in\operatorname{domain}\partial_{i}^{*}$
and that there exist $\zeta_{1},\dots,\zeta_{N}\in\mathcal{F}(R)$
so that \[
\zeta_{j}(X_{1},\dots,X_{N})=\partial_{j}(1\otimes1),\quad j=1,2,\dots,N.\]
Then there exists a II$_{1}$ factor $\mathcal{M}\cong M*L(\mathbb{F}_{\infty})$
and a $t_{0}>0$ so that for all $0\leq t<t_{0}$ there exists an
embedding $\alpha_{t}:M=W^{*}(X_{1},\ldots,X_{N})\to\mathcal{M}$
and a free $(0,1)$-semicircular family $S_{1},\ldots,S_{N}\in\mathcal{M}$,
free from $M$ and satisfying the inequality\begin{equation}
\Vert\alpha_{t}(X_{j})-(X_{j}+\sqrt{t}\Xi(X_{1},\dots,X_{N})\#S)\Vert_{2}\leq Ct,\label{eq:OWEstimate}\end{equation}
where $C$ is a fixed constant. Furthermore, $\alpha_{t}(X_{j})\in W^{*}(X_{1},\ldots,X_{N},S_{1},\ldots,S_{N},\{S_{j}'\}_{j=1}^{\infty})$,
where $\{S_{j}'\}_{j=1}^{\infty}$ are a free semicircular family,
free from $(X_{1},\ldots,X_{N},S_{1},\ldots,S_{N}).$ 

If $A$ can be embedded into $R^{\omega}$, so can $\mathcal{M}$.

In particular, the non-commutative Wasserstein distance of Biane-Voiculescu
satisfies:\[
d_{W}(\ (X_{j})_{j=1}^{N},\ (X_{j}+\sqrt{t}\Xi(X_{1},\dots,X_{N})\#S)_{j=1}^{N}\ )\leq Ct.\]

\end{prop}
\begin{proof}
By Lemma \ref{lem:1otimes1Enough}, we can find $\xi_{1},\dots,\xi_{N}\in\mathcal{F}(R)$
so that $\xi_{j}(X_{1},\dots,X_{N})=\partial^{*}\partial(X_{j})$. 

Let $\mathcal{M}=W^{*}(X_{1},\ldots,X_{N},\{S_{1}(s),\ldots,S_{N}(s):0\leq s\leq t\})$,
where $S_{j}(t)$ is a free semicircular Brownian motion. Let $X_{j}(t)$
be a stationary solution to the SDE \eqref{eq:SDE} (see Lemma \ref{lem:SDEstationary}).
The map that takes a polynomial in $X_{1},\ldots,X_{N}$ to a polynomial
in $X_{1}(t),\ldots,X_{N}(t)$ preserves traces and so extends to
an embedding $\alpha_{t}:M\to\mathcal{M}$. By the free Burkholder-Ghundy
inequality \cite{biane-speicher:stochcalc}, it follows that for $0\leq t<t_{0}<1$\[
\Vert X_{j}(t)-X_{j}(0)\Vert\leq C_{1}\sqrt{t}+C_{2}t\leq C_{3}\sqrt{t},\]
where $C_{1}=\sup_{t<t_{0}}\Vert\Xi(X_{1}(t),\dots,X_{N}(t)\Vert<\infty$,
$C_{2}=\max_{j}\sup_{t<t_{0}}\Vert\xi_{j}(X_{1},\dots,X_{N}(t)\Vert$$ $. 

Furthermore,\begin{eqnarray*}
X_{j}(t)-X_{j}(0) & = & \int_{0}^{t}\Xi(X_{1}(s),\ldots,X_{N}(s))\#dS_{j}(s)-\int_{0}^{t}\xi_{j}(X_{1}(s),\ldots,X_{Nn}(s))dt\\
 & = & \int_{0}^{t}\Xi(X_{1}(0),\ldots,X_{N}(0))\#dS_{j}(s)\\
 &  & -\int_{0}^{t}[\Xi(X_{1}(0),\ldots X_{N}(0))-\Xi(X_{1}(s),\ldots,X_{N}(s))]\#dS_{j}(s)\\
 &  & -\int_{0}^{t}\xi_{j}(X_{1}(s),\ldots,X_{N}(s))ds.\end{eqnarray*}
By the Lipschitz property of the coefficients of the SDE \eqref{eq:SDE},
we see that\[
\Vert\Xi(X_{1}(s),\ldots,X_{N}(s))-\Xi(X_{1}(0),\dots,X_{N}(0))\Vert\leq K\max_{j}\Vert X_{j}(s)-X_{j}(0)\Vert\leq K'\sqrt{s}.\]
Combining this with the estimate $\Vert\xi_{j}(X_{1}(t),\ldots,X_{N}(t))\Vert<K''$
we may apply the free Burkholder-Ghundy inequality to deduce that\begin{eqnarray*}
\Vert X_{j}(t)-(X_{j}(t)+\Xi(X_{1}(0),\ldots,X_{N}(0))\#S_{j}(t))\Vert & \leq & \left|\int_{0}^{t}(K'\sqrt{s})^{2}ds\right|^{1/2}+\Vert\int_{0}^{t}K''ds\Vert\\
 & \leq & Ct.\end{eqnarray*}
Thus it is enough to notice that $\Vert\cdot\Vert_{2}\leq\Vert\cdot\Vert$
and to take $S_{j}=\frac{1}{\sqrt{t}}S_{j}(t)$, which is a $(0,1)$
semicircular element.

If $M$ is $R^{\omega}$-embeddable, we may choose $\mathcal{M}$
to be $R^{\omega}$-embeddable as well, since it can be chosen to
be a free product of $M$ and a free group factor.

Finally, note that $X_{j}(t)\in W^{*}(X_{1},\ldots,X_{N},\{S_{j}(s):0\leq s\leq t\}_{j=1}^{N})$
by construction. But the algebra $W^{*}(\{S_{j}(s):0\leq s\leq t\})$
can be viewed as the algebra of the Free Gaussian functor applied
to the space $L^{2}[0,1]$, in such a way that $S_{j}(s)=S([0,s])$.
Then $W^{*}(\{S_{j}(s):0\leq s\leq t\})\subset W^{*}(S_{1},\ldots,S_{N},\{S_{k}'\}_{k\in I(j)})$,
where $\{S_{k}':k\in I(j)\}$ are free semicircular elements corresponding
to the completion of the singleton set $\{t^{-1/2}\chi_{[0,t]}\}$
to an ONB of $L^{2}[0,1]$.

The estimate for the Wasserstein distance now follows if we note that
the law of $(\alpha_{t}(X_{j}))_{j=1}^{N}$ is the same as that of
$(X_{j})_{j=1}^{N}$; thus $(X_{j}(t))_{j=1}^{N}\cup(X_{j}+\sqrt{t}\Xi\#S)_{j=1}^{N}$
is a particular $2N$-tuple with marginal distributions the same as
those of $(X_{j})_{j=1}^{N}$ and $(X_{j}+\sqrt{t}\Xi\#S)_{j=1}^{N}$,
so that the estimate \eqref{eq:OWEstimate} becomes an estimate on
the Wasserstein distance.
\end{proof}
\begin{rem}
Although we do not need this in the rest of the paper, we note that
the estimate in Proposition \ref{prop:OWprocess} also holds in the
operator norm.
\end{rem}
We should mention that an estimate similar to the one in Proposition
\ref{prop:OWprocess} was obtained by Biane and Voiculescu \cite{bine-voiculescu:WassersteinDist}
in the case $N=1$ under the much less restrictive assumptions that
$\Xi=1\otimes1$ and $1\otimes1\in\operatorname{domain}\partial^{*}$
(i.e., the free Fisher information $\Phi^{*}(X)$ is finite). Setting
$\Xi_{ij}=\delta_{ij}1\otimes1$ we have proved an analog of their
estimate (in the $N$-variable case), but under the very restrictive
assumption that the conjugate variables $\partial^{*}(\Xi)$ are analytic
functions in $X_{1},\ldots,X_{N}$. The main technical difficulty
in removing this restriction lies in the question of existence of
a stationary solution to \eqref{eq:SDE} in the case of very general
drifts $\xi$.

\section{Applications to $q$-semicircular families.}

\subsection{Estimates on certain operators related to $q$-semicircular families.}

\subsubsection{Background on $q$-semicircular elements.}

Let $H_{\mathbb{R}}$ be a finite-dimensional real Hilbert space,
$H$ its complexifiction $H=H_{\mathbb{R}}\otimes_{\mathbb{R}}\mathbb{C}$,
and let $F_{q}(H)$ be the $q$-deformed Fock space on $H$ \cite{Speicher:q-comm}.
Thus\[
F_{q}(H)=\mathbb{C}\Omega\oplus\bigoplus_{n\geq1}H^{\otimes n},\]
with the inner product given by\[
\langle\xi_{1}\otimes\cdots\otimes\xi_{n},\zeta_{1}\otimes\cdots\otimes\zeta_{m}\rangle=\delta_{n=m}\sum_{\pi\in S_{n}}q^{i(\pi)}\prod_{j=1}^{n}\langle\xi_{j},\zeta_{\pi(j))}\rangle,\]
where $i(\pi)=\#\{(i,j):i<j\ \textrm{and }\pi(i)>\pi(j)\}$.

We write $HS$ for the space of Hilbert-Schmidt operators on $F_{q}(H)$.
We denote by $\Xi\in HS$ the operator\[
\Xi=\sum q^{n}P_{n}\]
where $P_{n}$ is the orthogonal projection onto the subspace $H^{\otimes n}\subset F_{q}(H)$.

For $h\in H$, let $l(h):F_{q}(h)\to F_{q}(H)$ be the creation operator,
$l(h)(\xi_{1}\otimes\cdots\otimes\xi_{n})=h\otimes\xi_{1}\otimes\cdots\otimes\xi_{n}$,
and for $h\in H_{\mathbb{R}}$ let $s(h)=l(h)+l(h)^{*}$. We denote
by $M$ the von Neumann algebra $W^{*}(s(h):h\in H_{\mathbb{R}})$.
It is known \cite{ricard:qfactor,sniady:qfactor} that $M$ is a II$_{1}$
factor and that $\tau=\langle\cdot\Omega,\Omega\rangle$ is a faithful
tracial state on $M$. Moreover, $F_{q}(H)=L^{2}(M,\tau)$ and $HS=L^{2}(M,\tau)\bar{\otimes}L^{2}(M,\tau)$.

Fix an orthonormal basis $\{h_{i}\}_{i=1}^{N}\subset H_{\mathbb{R}}$
and let $X_{i}=s(h_{i})$. Thus $M=W^{*}(X_{1},\ldots,X_{N})$, $N=\dim H_{\mathbb{R}}$. 

\begin{lem}
\label{lem:oldResQ}\cite{shlyakht:qdim} For $j=1,\ldots,N$, let
$ $$ $$\partial_{j}:\mathbb{C}[X_{1},\ldots,X_{N}]\to HS$ be the
derivation given by $\partial_{j}(X_{i})=\delta_{i=j}\Xi$. Let $\partial:\mathbb{C}[X_{1},\ldots,X_{N}]\to HS^{N}$
be given by $\partial=\partial_{1}\oplus\cdots\oplus\partial_{N}$
and regard $\partial$ as an unbounded operator densely defined on
$L^{2}(M)$. Then:\\
(i) $\partial$ is closable. \\
(ii) If we denote by $Z_{j}$ the vector $0\oplus\cdots\oplus P_{\Omega}\oplus\cdots\oplus0\in HS^{N}$
(nonzero entry in $j$-th place, $P_{\Omega}$ is the orthogonal projection
onto $\mathbb{C}\Omega\in F_{q}(H)$), then $Z_{j}$ is in the domain
of $\partial^{*}$ and $\partial^{*}(Z_{j})=h_{j}$.
\end{lem}
As a consequence of (ii), if we let $\partial$ be as in the above
Lemma, $\xi_{j}=\partial^{*}(Z_{j})\in\mathbb{C}[X_{1},\cdots,X_{N}]\subset\mathcal{F}(R)$
for any $R$.

\subsubsection{$\Xi$ as an analytic function of $X_{1},\ldots,X_{n}$.}

We now claim that for small values of $q$, the element $\Xi\in L^{2}(M)^{\otimes2}$
defined in Lemma \ref{lem:oldResQ} can be thought of as an analytic
function of the variables $X_{1},\ldots,X_{N}$. Recall that $h_{i}\in H$
is a fixed orthonormal basis and $X_{j}=s(h_{j})$, $j=1,\dots,N$
thus form a $q$-semicircular family.

\begin{lem}
\label{lemma:Wick}Let $W_{i_{1},\ldots,i_{n}}$ be non-commutative
polynomials so that $W_{i_{1},\ldots,i_{n}}(X_{1},\ldots,X_{N})\Omega=h_{i_{1}}\otimes\cdots\otimes h_{i_{n}}$.
Then the degree of $W_{i_{1},\ldots,i_{n}}$ is $n$, and the maximal
absolute value $c_{k}^{(n)}$ of a coefficient of a monomial $X_{j_{1}}\cdots X_{j_{k}}$,
$k\leq n$, in $W_{i_{1},\ldots,i_{n}}$ satisfies\[
c_{k}^{(n)}\leq2^{n-k}\left(\frac{1}{1-|q|}\right)^{n-k}.\]
 Furthermore, $\Vert W_{i_{1},\ldots,i_{n}}\Vert_{L^{2}(M)}^{2}\leq2^{n}(1-|q|)^{-n}$.
\end{lem}
\begin{proof}
Clearly, $c_{n}^{(n)}=1$. Moreover, (compare \cite{effros-popa:qwick})\[
W_{i_{1},\ldots,i_{n}}=X_{i_{1}}W_{i_{2},\ldots,i_{n}}-\sum_{j\geq2}q^{j-2}\delta_{i_{1}=i_{j}}W_{i_{2},\ldots,\hat{i_{j}},\ldots,i_{n}}\]
(where $\hat{\cdot}$ denotes omission). So the degree of $W_{i_{1},\ldots,i_{n}}$
is $n$ and the coefficient $c_{n}$ of a monomial of degree $k$
in $W_{i_{1},\ldots,i_{n}}$ is at most the sum of a coefficient of
a degree $k-1$ monomial in $W_{i_{2},\ldots,i_{n}}$ and $\sum_{j\geq2}q^{j-2}|k_{j}|$,
where $k_{j}$ is a coefficient of a degree $k$ monomial in $W_{i_{2},\ldots,\hat{i}_{j},\ldots,i_{n}}$.
By induction, we see that\begin{eqnarray*}
c_{k}^{(n)} & \leq & c_{k-1}^{(n-1)}+\sum_{j\geq2}^{n}|q|^{j-2}c_{k}^{(n-2)}\\
 & \leq & 2^{n-k-2}\left(\frac{1}{1-|q|}\right)^{n-k}+2^{n-k-2}\left(\frac{1}{1-|q|}\right)^{n-k-2}\sum_{j\geq0}|q|^{j}\\
 & = & 2^{n-k-2}\left[\left(\frac{1}{1-|q|}\right)^{n-k}+\left(\frac{1}{1-|q|}\right)^{n-k-2}\frac{1}{1-|q|}\right]\\
 &  & \leq2^{n-k-2}\cdot2\left(\frac{1}{1-|q|}\right)^{n-k}\leq2^{n-k}\left(\frac{1}{1-|q|}\right)^{n-k}.\end{eqnarray*}
as claimed.

The upper estimate on $\Vert W_{i_{1},\ldots,i_{n}}\Vert_{L^{2}(M)}^{2}$
follows in a similar way.
\end{proof}
\begin{lem}
Let $\{\xi_{k}:k\in K\}$ be a finite set of vectors in an inner product
space $V$. Let $\Gamma$ be the matrix $\Gamma_{k,l}=\langle\xi_{k},\xi_{l}\rangle$.
Assume that $\Gamma$ is invertible and let $B=\Gamma^{-1/2}$. Then
the vectors\[
\zeta_{l}=\sum_{k}B_{k,l}\xi_{l}\]
form an orthonormal basis for the span of $\{\xi_{k}:k\in K\}$. Moreover,
if $\lambda$ denotes the smallest eigenvalue of $\Gamma$, then $|B_{k,l}|\leq\lambda^{-1/2}$
for each $k,l$.
\end{lem}
\begin{proof}
We have, using the fact that $B$ is symmetric and $B\Gamma B=I$:
$\langle\zeta_{l},\zeta_{l'}\rangle=\langle\sum_{k,k'}B_{k,l}\xi_{l},B_{k',l'}\xi_{l'}\rangle=\sum_{k,k'}B_{k,l}B_{k',l'}\Gamma_{l,l'}=(B\Gamma B)_{l,l'}=\delta_{l=l'}.$
\end{proof}
\begin{lem}
\label{lemma:ONB}There exist non-commutative polynomials $p_{i_{1},\ldots,i_{n}}$
in $X_{1},\ldots,X_{N}$ so that the vectors \[
\{p_{i_{1},\ldots,i_{n}}(X_{1},\ldots,X_{n})\Omega\}_{i_{1},\dots,i_{n}=1}^{N}\]
 are orthonormal and have the same span as $\{W_{i_{1},\ldots,i_{n}}\}_{i_{1},\ldots,i_{n}=1}^{N}$. 

Moreover, these can be chosen so that $p_{i_{1},\ldots,i_{n}}$ is
a polynomial of degree at most $n$ and the coefficient of each degree
$k$ monomial in $p$ is at most $(1-2|q|)^{-n/2}(2N)^{n}(1-|q|)^{k}2^{-k}$. 
\end{lem}
\begin{proof}
Consider the inner product matrix\[
\Gamma_{n}=[\langle W_{i_{1},\ldots,i_{n}},W_{j_{1},\ldots,j_{n}}\rangle]_{i_{1},\ldots,i_{n},\ j_{1},\ldots,j_{n}=1}^{N}.\]
Dykema and Nica proved (Lemma 3.1 \cite{dykema-nica:stability}) that
one has the following recursive formula for $\Gamma_{n}$. Consider
an $N^{n}$-dimensional vector space $W$ with orthonormal basis $e_{i_{1},\dots,i_{n}}$,
$i_{1},\dots,i_{n}\in\{1,\dots,N\}$, and consider the unitary representation
$\pi_{n}$ of the symmetric group $S_{n}$ given by $\sigma\cdot e_{i_{1},\dots,i_{n}}=e_{i_{\sigma(1)},\dots,i_{\sigma(n)}}$.
Denote by $(1\to j)$ the action (via $\pi_{n}$) of the permutation
that sends $1$ to $j$, $k$ to $k-1$ for $2\leq k\leq j$, and
$l$ to $l$ for $l>j$ on $W$. Let $M_{n}=\sum_{j=1}^{n}q^{j-1}(1\to j)\in\textrm{End}(W)$.
Then $\Gamma_{1}$ is the identity $N\times N$ matrix, and\[
\Gamma_{n}=(1\otimes\Gamma_{n-1})M_{n},\]
where $1\otimes\Gamma_{n}$ acts on the basis element $e_{j_{1},\dots,j_{n}}$
by sending it to $\sum_{k_{2},\dots,k_{n}}(\Gamma_{n-1})_{j_{2},\dots,j_{n},\ k_{2},\dots,k_{n}}e_{j_{1},k_{2},\dots,k_{n}}$
and $\Gamma$ acts on the basis elements by sending $e_{j_{1},\dots,j_{n}}$
to $\sum_{k_{1},\dots,k_{n}}(\Gamma_{n})_{j_{1},\dots,j_{n},\ k_{1},\dots,k_{n}}e_{k_{1},\dots,k_{n}}$.
They then proceeded to prove that the operator $M_{n}$ is invertible
and derive a bound for its inverse in the course of proving Lemma
4.1 in \cite{dykema-nica:stability}. Combining this bound and the
recursive formula for $\Gamma_{n}$ yields the following lower estimate
for the smallest eigenvalue of $\Gamma_{n}$:\begin{eqnarray*}
c_{n} & = & \left(\frac{1}{1-|q|}\prod_{k=1}^{\infty}\left(\frac{1-|q|^{k}}{1+|q|^{k}}\right)\right)^{n}=\left(\frac{1}{1-|q|}\sum_{k=-\infty}^{\infty}(-1)^{k}|q|^{k^{2}}\right)^{n}\\
 & \geq & \left(\frac{1}{1-|q|}\left[1-\sum_{k\geq0}|q|^{k^{2}}\right]\right)^{n}\geq\frac{1}{(1-|q|)^{n/2}}\left(1-\sum_{k\geq1}|q|^{k}\right)^{n}\\
 & \geq & \left(\frac{1}{1-|q|}\left[1-\frac{|q|}{1-|q|}\right]\right)^{n}=\left(\frac{1-2|q|}{\left(1-|q|\right)^{2}}\right)^{n}.\end{eqnarray*}

Thus if we set $B=\Gamma_{n}^{-1/2}$, then all entries of $B$ are
bounded from above by $c_{n}^{-1/2}$. Thus if we apply the previous
lemma with $K=\{1,\ldots,N\}^{n}$ to the vectors $\xi_{i_{1},\ldots,i_{n}}=W_{i_{1},\ldots,i_{n}}\Omega$,
we obtain that the vectors\[
\zeta_{i}=\sum_{j\in K}B_{j,i}\xi_{j},\quad i\in K\]
form an orthonormal basis for the subspace of the Fock space spanned
by tensors of length $n$.

Now for $i=(i_{1},\ldots,i_{n})\in K$, let\[
p_{i}(X_{1},\ldots,X_{N})=\sum_{j\in K}B_{j,i}W_{j}(X_{1},\ldots,X_{N}).\]
Then $\zeta_{i}=p_{i}(X_{1},\ldots,X_{N})\Omega$ are orthonormal
and (because the vacuum vector is separating), the polynomials $\{p_{i}:i\in K\}$
have the same span as $\{W_{i}:i\in K\}$.

Furthermore, if $a$ is the coefficient of a degree $k$ monomial
$r$ in $p_{i}$, then $a$ is a sum of at most $N^{n}$ terms, each
of the form $(\textrm{the coefficient of }r\mbox{ in }W_{j})B_{j,i}$.
Using Lemma \ref{lemma:Wick}, we therefore obtain the estimate\[
|a|\leq N^{n}c_{n}^{-1/2}2^{n-k}(1-|q|)^{-(n-k)}=\left(\frac{2N}{(1-2|q|)^{1/2}}\right)^{n}2^{-k}(1-|q|)^{k}.\]

\end{proof}
We'll now use the terminology of \S\ref{sub:PowerSeriesEst} in dealing
with non-commutative power series.

Let $R_{0}=2(1-|q|)^{-1}\geq2(1-q)^{-1}\geq\Vert X_{j}\Vert$. Then
if $\alpha>1$, $p=p_{i_{1},\dots,i_{n}}$ is as in Lemma \ref{lemma:ONB},
and $\phi_{p}$ is as in \S\ref{sub:PowerSeriesEst}, then the coefficient
of $z^{k}$, $k\leq n$ in $\phi_{p}$ is bounded by\[
\left(\frac{2N}{(1-2|q|)^{1/2}}\right)^{n}R_{0}^{-k}\leq\left(\frac{2N\alpha}{(1-2|q|)^{1/2}}\right)^{n}(\alpha NR_{0})^{-k}.\]
In particular for any $\rho<\alpha R_{0}$,\[
\Vert p_{i_{1},\dots,i_{n}}\Vert_{\rho}\le\left(\frac{2N\alpha}{(1-2|q|)^{1/2}}\right)^{n}\sum_{k=0}^{n}(\alpha NR_{0})^{-k}N^{k}\rho^{k}\leq\left(\frac{2N\alpha}{(1-2|q|)^{1/2}}\right)^{n}\frac{1}{1-\rho/(\alpha R_{0})}.\]

\begin{lem}
\label{lemma:XIisanalytic}Let $q$ be so that $|q|<(4N^{3}+2)^{-1}$.
Then:

(a) The formula\[
\Xi(Y_{1},\ldots,Y_{N})=\sum_{n}q^{n}\sum_{i_{1},\ldots,i_{n}}p_{i_{1},\ldots,i_{n}}(Y_{1},\ldots,Y_{N})\otimes p_{i_{1},\ldots,i_{n}}(Y_{1},\ldots,Y_{N})\]
defines a non-commutative power series with values in $\mathbb{C}\langle Y_{1},\dots,Y_{N}\rangle^{\otimes2}$
with radius of convergence strictly bigger than the norm of a $q$-semicircular
element, $\Vert X_{j}\Vert\leq2(1-q)^{-1}$.\\
(b) If $X_{1},\dots,X_{N}$ are $q$-semicircular elements and $\Xi$
is as in Lemma \ref{lem:oldResQ}, then $\Xi=\Xi(X_{1},\ldots,X_{N})$
(convergence in Hilbert-Schmidt norm, identifying $HS$ with $L^{2}(M)\bar{\otimes}L^{2}(M)$).
\end{lem}
\begin{proof}
Clearly,\[
\Vert p_{i_{1},\dots,i_{n}}\otimes p_{i_{1},\dots,i_{n}}\Vert_{\rho}\leq\Vert p_{i_{1},\dots,i_{n}}\Vert_{\rho}^{2}\leq\left(\frac{2N\alpha}{(1-2|q|)^{1/2}}\right)^{2n}\frac{1}{(1-\rho/(\alpha R_{0}))^{2}}=K_{\rho}\left(\frac{4N^{2}\alpha^{2}}{1-2|q|}\right)^{n}\]
for any $\rho<\alpha R_{0}$, where $R_{0}=2(1-|q|)^{-1}\geq\Vert X_{j}\Vert$. 

Thus\begin{eqnarray*}
\Vert\Xi\Vert_{\rho} & \leq & K_{\rho}\sum_{n}\left(\frac{4N^{2}\alpha^{2}}{1-2|q|}\right)^{n}|q|^{n}N^{n}\\
 & \leq & K_{\rho}\sum_{n}\left(\frac{4N^{3}\alpha|q|}{1-2|q|}\right)^{n},\end{eqnarray*}
which is finite as long as $\rho<\alpha R_{0}$ and \[
\frac{4N^{3}\alpha|q|}{1-2|q|}<1.\]
Thus as long as $4N^{3}|q|<1-2|q|$, i.e., $|q|<(4N^{3}+2)^{-1}$,
we can choose some $\alpha>1$ so that the series defining $\Xi$
has a radius of convergence of at least $\alpha R_{0}>\Vert X_{j}\Vert$. 

For part (b), we note that because $\Vert\cdot\Vert_{L^{2}(M)}\leq\Vert\cdot\Vert_{M}$
and because of the definition of the projective tensor product, we
see that\[
\Vert\cdot\Vert_{HS}\leq\Vert\cdot\Vert_{M\hat{\otimes}M}\]
on $M\hat{\otimes}M$. Thus convergence in the projective norm on
$M\hat{\otimes}M$ guarantees convergence in Hilbert-Schmidt norm.
Furthermore, by definition of orthogonal projection onto a space,\[
\Xi=\sum q^{n}P_{n}\]
where $P_{n}=\sum_{i_{1},\ldots,i_{n}}p_{i_{1},\ldots,i_{n}}\otimes p_{i_{1},\ldots,i_{n}}=\Xi^{(n)}(X_{1},\ldots,X_{N})$
are the partial sums of $\Xi(X_{1},\dots,X_{N})$ (here we again identify
$HS$ and $L^{2}\bar{\otimes}L^{2}$). Hence $\Xi=\Xi(X_{1},\ldots,X_{N})$.
\end{proof}

\section{An estimate on free entropy dimension.\label{sec:deltaEst}}

We now show how an estimate of the form \eqref{eq:OWlikeEstimate0}
can be used to prove a lower bound for the free entropy dimension
$\delta_{0}$. 

Recall \cite{dvv:entropy3,dvv:entropy2} that if $X_{1},\ldots,X_{n}\in(M,\tau)$
is an $n$-tuple of self-adjoint elements, then the set of microstates
$\Gamma_{R}(X_{1},\ldots,X_{n};l,k,\varepsilon)$ is defined by:\begin{multline*}
\Gamma_{R}(X_{1},\ldots,X_{n};l,k,\varepsilon)=\big{\{}(x_{1},\ldots,x_{n}\in(M_{k\times k}^{sa})^{n}:\Vert x_{j}\Vert<R,\\
\left|\tau(p(X_{1},\ldots,X_{n}))-\frac{1}{k}Tr(p(x_{1},\ldots,x_{n}))\right|<\varepsilon,\\
\textrm{for any monomial }p\mbox{ of degree }\leq l\big{\}}.\end{multline*}
If $R$ is omitted, the value $R=\infty$ is understood. 

The set of microstates for $X_{1},\ldots,X_{n}$ in the presence of
$Y_{1},\ldots,Y_{m}$ is defined by\begin{multline*}
\Gamma_{R}(X_{1},\ldots,X_{n}:Y_{1},\ldots,Y_{m};l,k,\varepsilon)=\big{\{}(x_{1},\ldots,x_{n}):\exists(y_{1},\ldots,y_{m})\\
\textrm{ s.t. }(x_{1},\ldots,x_{n},y_{1},\ldots,y_{m})\in\Gamma_{R}(X_{1},\dots,X_{n},Y_{1},\dots,Y_{m};l,k,\varepsilon)\big{\}}.\end{multline*}
The free entropy and free entropy in the presence are then defined
by \begin{eqnarray*}
\chi(X_{1},\ldots,X_{n}) & = & \sup_{R}\inf_{l,\varepsilon}\limsup_{k\to\infty}\left[\frac{1}{k^{2}}\log\textrm{Vol}(\Gamma_{R}(X_{1},\ldots,X_{n};l,k,\varepsilon))+\frac{n}{2}\log k\right]\\
\chi(X_{1},\ldots,X_{n}:Y_{1},\ldots,Y_{m}) & = & \sup_{R}\inf_{l,\varepsilon}\limsup_{k\to\infty}\Big[\frac{1}{k^{2}}\log\textrm{Vol}(\Gamma_{R}(X_{1},\ldots,X_{n}:Y_{1},\ldots,Y_{m};l,k,\varepsilon))\\
 &  & \qquad\qquad+\frac{n}{2}\log k\Big].\end{eqnarray*}
It is known \cite{bercovici-belinschi} that $\sup_{R}$ is attained;
in fact, $\chi(X_{1},\dots,X_{n}:Y_{1},\dots,Y_{m})=\chi_{R}(X_{1},\dots,X_{n}:Y_{1},\dots,Y_{m})$
once $R>\max_{i,j}\{\Vert X_{i}\Vert,\Vert Y_{j}\Vert\}$. 

Finally, the free entropy dimension $\delta_{0}$ is defined by\[
\delta_{0}(X_{1},\ldots,X_{n})=n+\limsup_{t\to0}\frac{\chi(X_{1}+\sqrt{t}S_{1},\ldots,X_{n}+\sqrt{t}S_{n}:S_{1},\ldots,S_{n})}{|\log t|},\]
where $S_{1},\ldots,S_{n}$ are a free semicircular family, free from
$X_{1},\ldots,X_{n}$. Equivalently \cite{jung:packing} one sets\[
\mathbb{K}_{\delta}(X_{1},\ldots,X_{n})=\inf_{\varepsilon,l}\limsup_{k\to\infty}\frac{1}{k^{2}}\log K_{\delta}(\Gamma_{\infty}(X_{1},\ldots,X_{n};k,l,\varepsilon)),\]
where $K_{\delta}(X)$ is the covering number of a set $X$ (the minimal
number of $\delta$-balls needed to cover $X$). Then\[
\delta_{0}(X_{1},\ldots,X_{n})=\limsup_{t\to0}\frac{\mathbb{K}_{t}(X_{1},\ldots,X_{n})}{|\log t|}.\]

\begin{lem}
\label{lem:QsharpScovering}Assume that $X_{1},\ldots,X_{n}\in(M,\tau)$,
$T_{jk}\in W^{*}(X_{1},\ldots,X_{n})\bar{\otimes}W^{*}(X_{1},\ldots,X_{n})^{op}$
are given. Set $S_{j}^{T}=\sum_{k}T_{jk}\#S_{k}$. Let $\eta=\dim_{M\bar{\otimes}M^{o}}(\overline{MS_{1}^{T}M+\cdots+MS_{n}^{T}M}^{L^{2}(M\bar{\otimes}M^{o})}).$

Then there exists a constant $K$ depending only on $T$ so that for
all $R>0$, $\alpha>0$, $t>0$, there are $\varepsilon'>0$, $l'>0$
and $k'>0$ so that for all $0<\varepsilon<\varepsilon'$, $k>k'$
and $l>l'$, and any $(x_{1},\ldots,x_{n})\in\Gamma(X_{1},\ldots,X_{n};k,l,\varepsilon)$
the set\begin{multline*}
\Gamma_{R}(tS_{1}^{I-T},\ldots,tS_{n}^{I-T}|(x_{1},\ldots,x_{n}):S_{1},\ldots,S_{n};k,l,\varepsilon)=\\
\{(y_{1},\ldots,y_{n}):\exists(s_{1},\ldots,s_{n})\textrm{ s.t. }(y_{1},\ldots,y_{n},x_{1},\ldots,x_{n},s_{1},\ldots,s_{n})\in\\
\Gamma_{R}(tS_{1}^{I-T},\ldots,tS_{n}^{I-T},X_{1},\ldots,X_{n},S_{1},\ldots,S_{n};k,l,\varepsilon)\}\end{multline*}
can be covered by $\left(K/t\right)^{(n-\eta+\alpha)k^{2}}$ balls
of radius $t^{2}$. 
\end{lem}
\begin{proof}
By considering the diffeomorphism of $(M_{k\times k}^{sa})^{n}$ given
by $(a_{1},\ldots,a_{n})\mapsto((1/t)a_{1},\ldots,(1/t)a_{n})$, we
may reduce the statement to showing that the set\[
\Gamma_{R}(S_{1}^{I-T},\ldots,S_{n}^{I-T}|(x_{1},\ldots,x_{n}):S_{1},\ldots S_{n};k,l,\varepsilon)\]
can be covered by $(C/t)^{(n-\eta+\alpha)k^{2}}$ balls of radius
$t$. 

Note that $\eta$ is the Murray-von Neumann dimension over $M\bar{\otimes}M^{o}$
of the image of the map $(\zeta_{1},\ldots,\zeta_{n})\mapsto(\zeta_{1}^{T},\ldots,\zeta_{n}^{T})$,
where $\zeta_{j}\in L^{2}(M)\bar{\otimes}L^{2}(M)$, $M=W^{*}(X_{1},\dots,X_{n})$.
Thus if we view $T$ as a matrix in $M_{n\times n}(M\bar{\otimes}M^{o})$,
then $\tau\otimes\tau\otimes Tr(E_{\{0\}}((I-T)^{*}(I-T)))=\eta$
(here $E_{X}$ denotes the spectral projection corresponding to the
set $X\subset\mathbb{R}$). 

Fix $\alpha>0$.

Then there exists $Q\in M_{n\times n}(\mathbb{C}[X_{1},\ldots,X_{n}]^{\otimes2})$
depending only on $t$ so that $\Vert Q_{ij}-(I-T)_{ij}\Vert_{2}<t/(2n)$. 

Set $S_{j}^{Q}=\sum_{k}Q_{jk}\#S_{k}$. Then\[
\Vert S_{j}^{Q(X_{1},\ldots,X_{n})}-S_{j}^{I-T}\Vert<\frac{t}{2}.\]
In particular, $\Vert S_{j}^{Q(X_{1},\ldots,X_{n})}-S_{j}^{I-T}\Vert_{2}<t/2$.
We may moreover choose $Q$ (again, depending only on $t$) so that\[
\tau\otimes\tau\otimes Tr(E_{[0,t/2[}(Q^{*}Q)^{1/2}(X_{1},\ldots,X_{n}))\geq\tau\otimes\tau\otimes Tr(E_{\{0\}}(I-T)^{*}(I-T))=\eta-\frac{1}{2}\alpha.\]

Thus for $l$ sufficiently large and $\varepsilon>0$ sufficiently
small, we have that if\[
(y_{1},\ldots,y_{n})\in\Gamma_{R}(S_{1}^{I-T},\ldots,S_{n}^{I-T}|(x_{1},\ldots,x_{n}):S_{1},\ldots,S_{n};k,l,\varepsilon),\]
then $\exists(s_{1},\ldots,s_{n})$ so that\[
(s_{1},\dots,s_{n},x_{1},\dots,x_{n})\in\Gamma_{R}(S_{1},\ldots,S_{n},X_{1},\ldots,S_{n};k,l,\varepsilon)\]
and\[
\Vert s_{j}^{Q(x_{1},\ldots,x_{n})}-y_{j}\Vert_{2}<t.\]
 By approximating the characteristic function $\chi_{[0,t/2]}$ with
polynomials on the interval $[0,\Vert Q(x_{1},\dots,x_{n})\Vert]$
(which is compact, since $\Vert x_{j}\Vert<R$), we may moreover assume
that $l$ is large enough and $\varepsilon$ is small enough so that\[
\frac{1}{k^{2}}Tr\otimes Tr\otimes Tr(E_{[0,t/2]}(Q^{*}Q)^{1/2}(x_{1},\ldots,x_{n}))\geq\eta-\alpha.\]
Denote by $\phi$ the map\[
(s_{1},\ldots,s_{n})\mapsto(s_{1}^{Q(x_{1},\ldots,x_{n})},\ldots,s_{n}^{Q(x_{1},\ldots,x_{n})}).\]
Let $R_{1}=\max_{j}\Vert S_{j}^{I-T}\Vert_{2}+1$. Assume that $\varepsilon<1$.
Then $\phi:(M_{k\times k}^{sa})^{n}\to(M_{k\times k}^{sa})^{n}$ is
a linear map, and since $\Vert s_{j}\Vert_{2}^{2}\leq1+\varepsilon<2$,
we have the inclusion:\[
\Gamma_{R}(S_{1}^{I-T},\ldots,S_{n}^{I-T}|(x_{1},\ldots,x_{n}):S_{1},\ldots,S_{n};k,l,\varepsilon)\subset N_{t}(\phi(B(2))\ \cap\ B(R_{1})),\]
where $B(R)$ the a ball of radius $R$ in $(M_{k\times k}^{sa})^{n}$
(endowed with the $L^{2}$ norm) and $N_{t}$ denotes a $t$-neighborhood.

The matrix of $\phi$ is precisely the matrix $Q(x_{1},\ldots,x_{n})\in M_{n\times n}(M_{k\times k})$. 

Let $\beta$ be such that $\beta nk^{2}$ eigenvalues of $(\phi^{*}\phi)^{1/2}$
are less than $R_{0}$. Then the $t$-covering number of $\phi(B(2))\cap B(R_{1})$
is at most\[
\left[\frac{R_{1}}{t}\right]^{(1-\beta)nk^{2}}\cdot\left[\frac{2R_{0}}{t}\right]^{\beta nk^{2}}.\]
Let $R_{0}=t/2$, so that $\beta=(\eta-\alpha)/n$. We conclude that
the $t$-covering number of $\Gamma_{R}(S_{1}^{I-T},\ldots,S_{n}^{I-T}|(x_{1},\ldots,x_{n}):S_{1},\ldots,S_{n};k,l,\varepsilon)$
is at most $(K/t)^{(n-\eta+\alpha)k^{2}}$, for some constant $K$
depending only on $R_{1}$, which itself depends only on $T$.
\end{proof}
\begin{thm}
\label{thm:deltaEstimate}Assume that $X_{1},\ldots,X_{n}\in(M,\tau)$,
$S_{1},\ldots,S_{n},\{S_{j}:j\in J\}$ is a free semicircular family,
free from $M$, $T_{jk}\in W^{*}(X_{1},\ldots,X_{n})\bar{\otimes}W^{*}(X_{1},\ldots,X_{n})^{op}$
are given, and that for each $t>0$ there exist $Y_{j}^{(t)}\in W^{*}(X_{1},\ldots,X_{n},S_{1},\ldots,S_{n},\{S_{j}'\}_{j\in J})$
so that:
\begin{itemize}
\item the joint law of $(Y_{1}^{(t)},\ldots,Y_{n}^{(t)})$ is the same as
that of $(X_{1},\ldots,X_{n})$;
\item If we set $S_{j}^{T}=\sum_{k}T_{jk}\#S_{k}$ and $Z_{j}^{(t)}=X_{j}+tS_{j}^{T}$,
then for some $t_{0}>0$ and some constant $C<\infty$ independent
of $t$, we have $\Vert Z_{j}^{(t)}-Y_{j}^{(t)}\Vert_{2}\leq Ct^{2}$
for all $t<t_{0}$.
\end{itemize}
Let $M=W^{*}(X_{1},\ldots,X_{n})$ and let\[
\eta=\dim_{M\bar{\otimes}M^{o}}(\overline{MS_{1}^{T}M+\cdots+MS_{n}^{T}M}^{L^{2}}).\]
Assume finally that $W^{*}(X_{1},\ldots,X_{n})$ embeds into $R^{\omega}$.
Then $\delta_{0}(X_{1},\ldots,X_{n})\geq\eta$.
\end{thm}
\begin{proof}
Let $T:(M\bar{\otimes}M^{o})^{n}\to(M\bar{\otimes}M^{o})^{n}$ be
the linear map given by\[
T(Y_{1},\ldots,Y_{n})=(\sum_{k}T_{1k}\#Y_{k},\ldots,\sum_{k}T_{nk}\#Y_{k}).\]
Then $\eta$ is the Murray-von Neumann dimension of the image of $T$,
and consequently\[
\eta=n-\dim_{M\bar{\otimes}M^{o}}\ker T.\]

Let $t$ be fixed. 

Since $Y_{j}^{(t)}$ can be approximated by non-commutative polynomials
in $X_{1},\ldots,X_{n}$, $S_{1},\ldots,S_{n}$ and $\{S_{j}':j\in J\}$,
for any $k_{0},\varepsilon_{0},l_{0}$ sufficiently large we may find
$k>k_{0}$, $l>l_{0}$, $\varepsilon<\varepsilon_{0}$ and $J_{0}\subset J$
finite so that whenever \[
(z_{1},\ldots,z_{n})\in\Gamma_{R}(X_{1}+tS_{1}^{T},\ldots,X_{n}+tS_{n}^{T}:S_{1},\ldots,S_{n},\{S_{j}'\}_{j\in J_{0}};k,l,\varepsilon),\]
there exists\[
(y_{1},\ldots,y_{n})\in\Gamma_{R}(X_{1},\ldots,X_{n};k,l_{0},\varepsilon_{0})\]
so that\begin{eqnarray}
\Vert y_{j}-z_{j}\Vert_{2} & \leq & Ct^{2}.\label{eq:yzOW}\end{eqnarray}

For a set $X\subset(M_{k\times k}^{sa})^{n}$ we'll write $K_{r}$
for its covering number by balls of radius $r$.

Consider a covering of $\Gamma_{R}(X_{1}+tS_{1},\ldots,X_{n}+tS_{n}:S,\ldots,S_{n},\{S_{j}'\}_{j\in J_{0}};k,l,\varepsilon)$
by balls of radius $(C+2)t^{2}$ constructed as follows.

First, let $(B_{\alpha})_{\alpha\in I}$ be a covering of $\Gamma_{R}(X_{1}+tS_{1}^{T},\ldots,X_{n}+tS_{n}^{T}:S_{1},\ldots,S_{n},\{S_{j}'\}_{j\in J_{0}};k,l_{0},\varepsilon_{0})$
by balls of radius $(C+1)t^{2}$. Because of \eqref{eq:yzOW}, we
may assume that\[
|I|\leq K_{t^{2}}(\Gamma_{R}(X_{1},\ldots,X_{n};k,l,\varepsilon)).\]
Next, for each $\alpha\in I$, let $(x_{1}^{(\alpha)},\ldots,x_{n}^{(\alpha)})\in B_{\alpha}$
be the center of $B_{\alpha}$. Consider a covering $(C_{\beta}^{(\alpha)}:\beta\in J_{\alpha})$
of $\Gamma_{R}(tS_{1}^{I-T},\ldots,tS_{n}^{I-T}|(x_{1}^{(\alpha)},\ldots,x_{n}^{(\alpha)}):S_{1},\ldots,S_{n},)$
by balls of radius $t^{2}$. By Lemma \ref{lem:QsharpScovering},
this covering can be chosen to contain $|J_{\alpha}|\leq(K/t)^{n-\eta'}$
balls, for any $\eta'<\eta$. Thus the sets\[
(B_{\alpha}+C_{\beta}^{(\alpha)}:\alpha\in I,\beta\in J_{\alpha}),\]
each of which is contained in a ball of radius at most $(C+2)t^{2}$,
cover the set $\Gamma_{R}(X_{1}+tS_{1},\ldots,X_{n}+tS_{n}:S_{1},\ldots,S_{n};k,l_{0},\varepsilon_{0})$.
The cardinality of this covering is at most\[
f(t^{2},k)\leq|I|\cdot\sup_{\alpha}|J_{\alpha}|\leq K_{t^{2}}(\Gamma_{R}(X_{1},\ldots,X_{n};k,l,\varepsilon)\cdot(Kt)^{\eta'-n}.\]

It follows that if we denote by $V(R,d)$ the volume of a ball of
radius $R$ in $\mathbb{R}^{d}$, we find that\[
\operatorname{Vol}(\Gamma_{R}(X_{1}+tS_{1},\ldots,X_{n}+tS_{n}:S_{1},\ldots,S_{n},\{S_{j}'\}_{j\in J_{0}}))\leq f(t^{2},k)\cdot V((C+2)t^{2},nk^{2}),\]
so that if we denote by $\mathbb{K}_{t^{2}}(X_{1},\ldots,X_{n})$
the expression $\inf_{\varepsilon,l}\limsup_{k\to\infty}\frac{1}{k^{2}}\log K_{t^{2}}(\Gamma(X_{1},\ldots,X_{n};k,l,\varepsilon))$
and set $C'=\log(C+2)$, we obtain the inequality\begin{multline*}
\inf_{\epsilon,l}\limsup_{k\to\infty}\frac{1}{k^{2}}\log\textrm{Vol}\Gamma_{R}(X_{1}+tS_{1},\ldots,X_{n}+tS_{n}:S_{1},\ldots,S_{n},\{S_{j}'\}_{j\in J_{0}};k,l,\varepsilon)\\
\leq\limsup_{k\to\infty}\log f(t^{2},k)+2n\log t+\log(C+2)\\
\leq\mathbb{K}_{t^{2}}(X_{1},\ldots,X_{n})+(\eta'-n)\log Kt+2n\log t+C'\\
=\mathbb{K}_{t^{2}}(X_{1},\ldots,X_{n})+(\eta'+n)\log t+(\eta'-n)\log K+C'.\end{multline*}
By freeness of $\{S_{j}'\}_{j\in J}$ and $\{S_{1},\ldots,S_{n},X_{1},\ldots,X_{n}\}$,
the $\limsup$ on the right-hand side remains the same if we take
$J_{0}=\emptyset$. Thus\[
\chi_{R}(X_{1}+tS_{1},\ldots,X_{n}+tS_{n}:S_{1},\dots,S_{n})\leq\mathbb{K}_{t^{2}}(X_{1},\ldots,X_{n})+(\eta'+n)\log t+C''.\]
If we divide both sides by $|\log t|$ and add $n$ to both sides
of the resulting inequality, we obtain:\begin{eqnarray*}
n+\frac{\chi_{R}(X_{1}+tS_{1},\ldots,X_{n}+tS_{n}:S_{1},\ldots,S_{n})}{|\log t|} & \leq & \frac{\mathbb{K}_{t^{2}}(X_{1},\ldots,X_{n})}{|\log t|}+(\eta'+n)\frac{\log t}{|\log t|}+n\\
 & = & 2\frac{\mathbb{K}_{t^{2}}(X_{1},\ldots,X_{n})}{|\log t^{2}|}+(\eta'+n)\frac{\log t}{|\log t|}+n.\end{eqnarray*}
Taking $\sup_{R}$ and $\limsup_{t\to0}$ and noticing that $\log t<0$
for $t<1$, we get the inequality\[
\delta_{0}(X_{1},\ldots,X_{n})\leq2\delta_{0}(X_{1},\ldots,X_{n})-(\eta+n)+n=2\delta_{0}(X_{1},\ldots,X_{n})-\eta'.\]
Solving this inequality for $\delta_{0}(X_{1},\ldots,X_{n})$ gives
finally\[
\delta_{0}(X_{1},\ldots,X_{n})\geq\eta'.\]
Since $\eta'<\eta$ was arbitrary, we obtain that $\delta_{0}(X_{1},\ldots,X_{n})\geq\eta$
as claimed.
\end{proof}
\begin{cor}
Let $(A,\tau)$ be a finitely-generated algebra with a positive trace
$\tau$ and generators $X_{1},\ldots,X_{n}$, and let $\operatorname{Der}_{a}(A;A\otimes A)$
denote the space of derivations from $A$ to $L^{2}(A\otimes A,\tau\otimes\tau)$
which are $L^{2}$ closable and so that for some $\Xi_{j}\in\mathcal{F}'(R)$,
$\xi\in\mathcal{F}(R)$, $R>\max_{j}\Vert X_{j}\Vert$, $\partial^{*}(1\otimes1)=\xi(X_{1},\dots,X_{n})$
and $\partial(X_{j})=\Xi_{j}(X_{1},\dots,X_{n})$. Consider the A,A-bimodule\[
V=\{(\delta(X_{1}),\ldots,\delta(X_{n})):\delta\in\operatorname{Der}_{a}(A;A\otimes A)\}\subset L^{2}(A\otimes A,\tau\otimes\tau)^{n}.\]
Assume finally that $M=W^{*}(A,\tau)\subset R^{\omega}$. Then\[
\delta_{0}(X_{1},\ldots,X_{n})\geq\dim_{M\bar{\otimes}M^{o}}\overline{V}^{L^{2}(A\otimes A,\tau\otimes\tau)^{n}}.\]

\end{cor}
\begin{proof}
Let $P:L^{2}(A\otimes A,\tau\otimes\tau)^{n}\to\overline{V}$ be the
orthogonal projection, and let $v_{j}=P(0,\ldots,1\otimes1,\ldots,0)$
($1\otimes1$ in the $j$-th position). Let $v_{j}^{(k)}=(v_{1j}^{(k)},\ldots,v_{nj}^{(k)})\in L^{2}(A\otimes A)^{n}$
be vectors approximating $v_{j}$, having the property that the derivations
defined by $\delta(X_{j})=v_{ij}^{(k)}$ lie in $\operatorname{Der}_{a}$.
Then\[
\eta_{k}=\dim_{M\bar{\otimes}M^{o}}\overline{Av_{1}^{(k)}A+\cdots+Av_{n}^{(k)}A}\to\dim_{M\bar{\otimes}M^{o}}\overline{V}\]
as $k\to\infty$. Now for each $k$, the derivations $\delta_{j}:A\to L^{2}(A\otimes A)$
so that $\delta_{j}(X_{i})=v_{ij}^{(k)}$ belong to $\operatorname{Der}_{a}$.
Applying Lemma \ref{lem:1otimes1Enough} and Proposition \ref{prop:OWprocess}
to $T_{ij}=v_{ij}^{(k)}$ and combining the conclusion with Theorem
\ref{thm:deltaEstimate} gives that\[
\delta_{0}(X_{1},\ldots,X_{n})\geq\eta_{k}.\]
Taking $k\to\infty$ we get\[
\delta_{0}(X_{1},\ldots,X_{n})\geq\dim_{M\bar{\otimes}M^{o}}V,\]
as claimed.
\end{proof}
\begin{cor}
For a fixed $N$, and all $|q|<(4N^{3}+2)^{-1}$, the $q$-semicircular
family $X_{1},\ldots,X_{N}$ satisfies\[
\delta_{0}(X_{1},\ldots,X_{N})>1\textrm{ and }\delta_{0}(X_{1},\ldots,X_{N})\geq N\left(1-\frac{q^{2}N}{1-q^{2}N}\right).\]
In particular, $M=W^{*}(X_{1},\ldots,X_{N})$ has no Cartan subalgebra.
Moreover, for any abelian subalgebra $\mathcal{A}\subset M$, $L^{2}(M)$,
as an $\mathcal{A},\mathcal{A}$-bimodule, contains a copy of the
coarse correspondence.
\end{cor}
\begin{proof}
Let $\partial_{i}$ be a derivation as in Lemma \ref{lem:oldResQ};
thus $\partial_{i}(X_{j})=\delta_{i=j}\Xi$, as defined in Lemma \ref{lem:oldResQ}.
Then for $|q|<(4N^{3}+2)^{-1}$, Lemma \ref{lemma:XIisanalytic} shows
that $\partial_{i}\in\operatorname{Der}_{a}$. Then Theorem \ref{thm:deltaEstimate}
implies that\[
\delta_{0}(X_{1},\ldots,X_{N})\geq\dim_{M\bar{\otimes}M^{o}}\overline{\sum M\Xi_{i}M},\]
$M=W^{*}(X_{1},\ldots,X_{n})$. It is known \cite{shlyakht:qdim}
that for $q^{2}<1/N$ (which is the case if we make the assumptions
about $q$ as in the hypothesis of the corollary), this dimension
is strictly bigger than $1$, and is no less than $N(1-q^{2}N(1-q^{2}N)^{-1})$. 

The facts about $M$ follow from Voiculescu's results \cite{dvv:entropy3}.
\end{proof}
For $N=2$, $(4N^{3}+2)^{-1}=1/34$. Thus the theorem applies for
$0\leq q\leq1/34=0.029\ldots$. Our estimate also shows that as $q\to0$,
$\delta_{0}(X_{1},\ldots,X_{N})\to N$. 

\begin{cor}
Let $\Gamma$ be a discrete group generated by $g_{1},\dots,g_{n}$,
and let $V\subset C^{1}(\Gamma,\ell^{2}\Gamma)$ be the subset consisting
of cocycles valued in $\mathbb{C}\Gamma\subset\ell^{2}\Gamma$. If
the group von Neumann algebra of $\Gamma$ can be embedded into the
ultrapower of the hyperfinite II$_{1}$ factor (e.g., if the group
is soffic), then\[
\delta_{0}(\mathbb{C}\Gamma)\geq\dim_{L(\Gamma)}\overline{V}.\]

\end{cor}
\begin{proof}
Any such cocycle gives rise to a derivation into $\mathbb{C}\Gamma^{\otimes2}$
by the formula\[
\partial(\gamma)=c(\gamma)\otimes\gamma^{-1}.\]
Then $\partial^{*}\partial(\gamma)=\Vert c(\gamma)\Vert_{2}^{2}\gamma\in\mathbb{C}\Gamma$.
Moreover, the bimodule generated by values of these derivations on
any generators of $\mathbb{C}\Gamma$ has the same dimension over
$L(\Gamma)\bar{\otimes L(\Gamma)}$ as $\dim_{L(\Gamma)}\bar{V}$.
\end{proof}
For certain $R^{\omega}$ embeddable groups (e.g., free groups, amenable
groups, residually finite groups with property $T$, more generally
embeddable groups with first $L^{2}$ Betti number $\beta_{1}^{(2)}=0$,
as well as groups obtained from these by taking amalgamated free products
over finite subgroups and passing to finite index subgroups and finite
extensions), $V$ is actually dense in the set of $\ell^{2}$ 1-cocycles.
Indeed, this is the case if all $\ell^{2}$ derivations are inner
(i.e., $\beta_{1}^{(2)}(\Gamma)=0$). Moreover, it follows from the
Meyer-Vietoris exact sequence that amalgamated free products over
finite subgroups retain the property that $V$ is dense in the space
of $\ell^{2}$ cocycles. Moreover, this property is also clearly preserved
by passing to finite-index subgroups and finite extensions. So it
follows that for such groups $\Gamma$, $\delta_{0}(\Gamma)=\beta_{1}^{(2)}(\Gamma)+\beta_{0}^{(2)}(\Gamma)-1$
(compare \cite{brown-dykema-jung:fdimamalg}). 

It is likely that the techniques of the present paper could be extended
to answer in the affirmative the following:

\begin{conjecture}
Let $\Gamma$ be a group generated by $g_{1},\dots,g_{n}$ and assume
that $L(\Gamma)$ can be embedded into $R^{\omega}$. Let $V\subset\ell^{2}(\Gamma)^{n}$
be the subspace $\{(c(g_{1}),\dots,c(g_{n})):c:\Gamma\to\ell^{2}(\Gamma)\ 1\textrm{-cocycle}\}$.
Let $P_{V}:\ell^{2}(\Gamma)^{n}\to V$ be the orthogonal projection,
so that $P_{V}\in M_{n\times n}(R(\Gamma))$, where $R(\Gamma)$ is
the von Neumann algebra generated by the right regular representation
of the group.

Let $\mathcal{A}\subset R(\Gamma)$ be the closure of $\mathbb{C}\Gamma\subset R(\Gamma)$
under holomorphic functional calcuclus, and let $P_{a}\in\mathcal{A}$
be any projection so that $P_{a}\leq P_{V}$. Then $\delta_{0}(\Gamma)\geq Tr_{M_{n\times n}}\otimes\tau_{R(\Gamma)}(P_{a})$.
\end{conjecture}
Note that with the notations of the Conjecture, $Tr_{M_{n\times n}}\otimes\tau_{R(\Gamma)}(P_{V})=\beta_{1}^{(2)}(\Gamma)-\beta_{0}^{(2)}(\Gamma)+1=\delta^{*}(\Gamma)$.

It should be noted that the restriction on the values of the cocycles
($\mathbb{C}\Gamma$ rather than $\ell^{2}\Gamma$) comes from the
difficulty in the extending the results of Proposition \ref{prop:OWprocess}
to the case of non-analytic $\Xi$ (though the term $\partial^{*}\partial(\gamma)$
continues to be a polynomial even in the case that the cocycle is
valued in $\ell^{2}(\Gamma)$ rather than $\mathbb{C}\Gamma$).

\section{Appendix: Otto-Villani type estimates via exponentiation of derivations.\label{sec:owDer}}

Let $M=W^{*}(X_{1},\ldots,X_{N})$, where $X_{1},\ldots,X_{N}$ are
self-adjoint. 

Let us denote by $\zeta_{j}$ the vector $(0,\ldots,0,1\otimes1,0,\ldots,0)\in[L^{2}(M,\tau)^{\otimes2}]^{N}$
(the only non-zero entry is in the $j$-th position). One can realize
a free semicircular family of cardinality $N$ on the space\[
H=L^{2}(M,\tau)\oplus\bigoplus_{k\geq1}\left[(L^{2}(M,\tau)\bar{\otimes}L^{2}(M,\tau))^{\oplus N}\right]^{\otimes_{M}k}.\]
using creation and annihilation operators: $S_{i}=L_{i}+L_{i}^{*}$
where\[
L_{i}\xi=\zeta_{i}\otimes_{M}\xi.\]
Then for $\zeta\in W^{*}(M)\bar{\otimes}W^{*}(M)$, the notation $S_{\zeta}$
makes sense, with $S_{\zeta_{i}}=S_{i}$, $aS_{\zeta}b+b^{*}S_{\zeta}a^{*}=S_{a\zeta b+b^{*}\zeta a^{*}}$
and $\Vert S_{\zeta}\Vert_{2}=\Vert\zeta\Vert_{2}$. 

Let $A=\operatorname{Alg}(X_{1},\ldots,X_{N})$. 

For $a,b\in A\otimes A$ and $j=1,\ldots,N$ let us write\[
(a\otimes b)\#S=aSb.\]

With these notations, we have:

\begin{prop}
\label{prop:canExp}Let $\partial:A\to V_{0}=[W^{*}(M,\tau)\bar{\otimes}W^{*}(M,\tau)]^{\oplus N}\subset V=[L^{2}(M,\tau)\bar{\otimes}L^{2}(M,\tau)]^{\oplus N}$
be a derivation. We assume that for each $j$, $\zeta_{j}$ is in
the domain of $\partial^{*}:V\to L^{2}(M,\tau)$ and that $\partial(a^{*})=(\partial(a))^{*}$,
where $*:L^{2}(M)\bar{\otimes}L^{2}(M)$ is the involution $(a\otimes b)^{*}=b^{*}\otimes a^{*}$.
Let $S_{1},S_{2},\ldots$ be semicircular elements, free from $M$. 

Assume that $\partial(A)\subset(A\otimes A)^{\oplus N}$ and also
that $\partial^{*}(1\otimes1)\in A$.

Then there exists a one-parameter group $\alpha_{t}$ of automorphisms
of $M*W^{*}(S_{1},\dots,S_{N})\cong M*L(\mathbb{F}_{N})$ so that
$A\cup\{S_{j}:1\leq j\leq N\}$ are analytic for $\alpha_{t}$ and
\begin{eqnarray*}
\frac{d}{dt}\Big\vert_{t=0}\alpha_{t}(a) & = & S_{\partial(a)},\qquad\forall a\in A,\\
\frac{d}{dt}\Big\vert_{t=0}\alpha_{t}(S_{j}) & = & -\partial^{*}(\zeta_{j}),\qquad j=1,2,\ldots.\end{eqnarray*}
In particular,\[
\alpha_{t}(a)\cdot1=(a-\frac{t^{2}}{2}\sum_{j}\partial^{*}(\partial(a))\ +\ t\partial(a)\ -\ \frac{t^{2}}{2}(1\otimes\partial+\partial\otimes1)(\partial(a))\in H.\]

\end{prop}
\begin{proof}
Let $B$ be the algebra generated by $A$ and $S_{1},\ldots,S_{N}$
in $\mathcal{M}=W^{*}(A,\tau)*L(\mathbb{F}_{N})$. 

Let $P_{j}:V\to L^{2}(A\otimes A)$ be the $j$-th coordinate projection,
and let $\partial_{j}:A\to A\otimes A$ be given by $\partial_{j}=P_{j}\circ\partial$.$ $

Let $V_{1},\ldots,V_{N}\in B$ be given by\[
V_{j}=\sum_{k}\partial_{k}(X_{j})\#S_{k}=S_{\partial(X_{j})},\qquad j=1,\ldots,N.\]
Let $V_{N+1},\ldots,V_{2N}\in B$ be given by\[
V_{N+k}=-\partial_{k}^{*}(1\otimes1)=-\partial^{*}(\zeta_{k}),\qquad k=1,\ldots,N.\]
Then $(V_{1},\ldots,V_{2N})\in B\subset L^{2}(B,\tau)$ is a non-commutative
vector field in the sense of \cite{dvv:cyclomorphy}. It is routine
to check that this vector field is orthogonal to the cyclic gradient
space. 

We now use \cite{dvv:cyclomorphy} to deduce that there exists a one-parameter
automorphism group $\alpha_{t}$ of $\mathcal{M}=W^{*}(B,\tau)$ with
the property that\begin{eqnarray*}
\frac{d}{dt}\alpha_{t}(X_{j})\Big|_{t=0} & = & V_{j},\qquad j=1,\ldots,N\\
\frac{d}{dt}\alpha_{t}(S_{k})\Big|_{t=0} & = & V_{N+k},\qquad k=1,\ldots,N,\end{eqnarray*}
and moreover that all elements in $B$ are analytic for $\alpha_{t}$.
In particular, we see that\begin{eqnarray*}
\frac{d}{dt}\alpha_{t}(X_{j})\Big|_{t=0} & = & S_{\partial(X_{j})},\\
\frac{d^{2}}{dt^{2}}\alpha_{t}(X_{j})\Big|_{t=0} & \cdot1= & \delta(S_{\partial(X_{j})})=-\partial^{*}(\partial(X_{j}))-(1\otimes\partial+\partial\otimes1)(\partial(X_{j})),\end{eqnarray*}
as claimed.
\end{proof}
\begin{example}
We give three examples in which the automorphisms $\alpha_{t}$ can
be explicitly constructed. The first is the case that $X_{1},\ldots,X_{N}$
is a free semicircular system and $\partial(X_{j})=(0,\ldots,1\otimes1,\ldots0)$
(i.e., $\partial=\oplus\partial_{j}$, where $\partial_{j}$ are the
difference quotient derivations of \cite{dvv:entropy5}). In this
case, the automorphism $\alpha_{t}$ is given by\[
\alpha_{t}(X_{j})=(\cos t)X_{j}+(\sin t)S_{j},\qquad\alpha_{t}(S_{j})=-(\sin t)X_{j}+(\cos t)S_{j}.\]
Another situation is that of a general $N$-tuple $X_{1},\ldots,X_{N}$
and $\partial$ an inner derivation given by $\partial(X)=[X,T]$,
for $[T_{j}]_{j=1}^{N}=[-T_{j}^{*}]_{j=1}^{N}\in[M\bar{\otimes}M^{o}]^{N}$.
Put $z=\sum T_{j}\#S_{j}$. Then $\alpha_{t}$ is an inner automorphism
given by $\alpha_{t}(Y)=\exp(izt)Y\exp(-izt)$. Lastly, assume that
$M=M_{1}*M_{2}$ and the derivations $\partial_{j}$ are determined
by $\partial_{j}|_{M_{1}}=0$, $\partial_{j}|_{M_{2}}(x)=[x,T_{j}]$
for some $T_{j}\in M\bar{\otimes}M^{o}$. Then again put $z=\sum T_{j}\#S_{j}$.
The automorphism $\alpha_{t}$ is then given by $\alpha_{t}(Y)=\exp(izt)Y\exp(-izt)$.
In particular, $\alpha_{t}|_{M_{1}}=\textrm{id}$ and $\alpha_{t}|_{M_{2}}$
is given by conjugation by unitaries $\exp(izt)$ which are free from
$M_{1}$ and $M_{2}$.
\end{example}
Proposition \ref{prop:canExp} can be used to give another proof to
the Otto-Villani type estimates (Proposition \ref{prop:OWprocess})
in the case of polynomial coefficients, using the following standard
lemma:

\begin{lem}
\label{lem:autEst}Let $\beta_{t}:(M,\tau)\to(M,\tau)$ be a one-parameter
group of automorphisms so that $\tau\circ\beta_{t}=\tau$. Let $X\in M$
be an element so that $t\mapsto\beta_{t}(X)$ is twice-differentiable.
Finally let $Z=\frac{d}{dt}\beta_{t}(X)\Big|_{t=0}$, $\xi=\frac{d^{2}}{dt^{2}}\beta_{t}(X)\Big|_{t=0}$.

Then one has for all $t$\[
\Vert\beta_{t}(X)-(X+tZ)\Vert_{2}\leq\frac{t^{2}}{2}\Vert\xi\Vert_{2}.\]

\end{lem}
\begin{cor}
Assume that $X_{1},\ldots,X_{N}\in A$ and $\partial_{1},\ldots\partial_{N}:A\to A\otimes A$
are derivations, so that $\partial_{j}^{*}(1\otimes1)\in A$. Then
we have the following estimate for the free Wasserstein distance:\[
d_{W}((X_{1},\ldots,X_{N}),(X_{1}+\sqrt{t}\sum_{k}\partial_{k}(X_{1})\#S_{k},\ldots,X_{N}+\sqrt{t}\sum_{k}\partial_{k}(X_{N})\#S_{k}))\leq Ct\]
where $C$ is the constant given by\[
C=\frac{1}{2}\left(\sum_{j}\Vert\partial^{*}\partial(X_{j})\Vert_{L^{2}(A)}^{2}+\Vert(1\otimes\partial+\partial\otimes1)(\partial(X_{j}))\Vert_{[L^{2}(A)\otimes L^{2}(A)\otimes L^{2}(A)]^{N^{2}}}^{2}\right)^{1/2},\]
where $\partial:A\to[L^{2}(A)\otimes L^{2}(A)]^{N}$ is the derivation
$\partial=\partial_{1}\oplus\cdots\oplus\partial_{N}$.

In the specific case of the difference quotient derivations determined
by $\partial_{k}(X_{j})=\delta_{kj}1\otimes1$, we have\[
d_{W}((X_{1},\ldots,X_{N}),(X_{1}+\sqrt{t}S_{1},\ldots,X_{N}+\sqrt{t}S_{N}))\leq\frac{t}{2}\Phi^{*}(X_{1},\ldots,X_{N})^{1/2}.\]

\end{cor}
\begin{proof}
Let $\alpha_{t}$ be the one-parameter group of automorphisms as in
Proposition \ref{prop:canExp}. We note that\[
\left(\sum_{j}\Vert\alpha_{\sqrt{t}}(X_{j})-(X_{j}+\sqrt{t}\sum_{k}\partial_{k}(X_{j})\#S_{k})\Vert_{2}^{2}\right)^{1/2}\leq Ct\]
in view of Lemma \ref{lem:autEst} and the expression for $\alpha''_{t}(X_{j})$.
On the other hand, $(\alpha_{\sqrt{t}}(X_{1}),\ldots,\alpha_{\sqrt{t}}(X_{N}))$
has the same law as $(X_{1},\ldots,X_{N})$, since $\alpha_{\sqrt{t}}$
is a $*$-homomorphism. It therefore follows that\begin{gather*}
d_{W}(X_{1},\ldots,X_{N},(X_{1}+\sqrt{t}\sum_{k}\partial_{k}(X_{1})\#S_{k},\ldots,X_{N}+\sqrt{t}\sum_{k}\partial_{k}(X_{N})\#S_{k}))\\
=d_{W}(\alpha_{\sqrt{t}}(X_{1}),\ldots,\alpha_{\sqrt{t}}(X_{N}),(X_{1}+\sqrt{t}\sum_{k}\partial_{k}(X_{1})\#S_{k},\ldots,X_{N}+\sqrt{t}\sum_{k}\partial_{k}(X_{N})\#S_{k}))\\
\leq Ct.\end{gather*}
In the case of the difference quotient derivations, we have:\begin{gather*}
\sum_{k}\partial_{k}(X_{j})\#S_{k}=S_{j};\qquad(1\otimes\partial+\partial\otimes1)(\partial(X_{j}))=(1\otimes\partial+\partial\otimes1)(1\otimes1)=0;\\
\partial^{*}\partial(X_{j})=\partial_{j}^{*}(1\otimes1).\end{gather*}
Thus\[
C=\frac{1}{2}\left(\sum_{j}\Vert\partial_{j}^{*}(1\otimes1)\Vert_{2}^{2}\right)^{1/2}=\frac{1}{2}\Phi^{*}(X_{1},\ldots,X_{N})^{1/2}\]
as claimed.
\end{proof}
\bibliographystyle{amsplain}
%\bibliography{/Volumes/shlyakht/tex/bib/quasifree}

\providecommand{\bysame}{\leavevmode\hbox to3em{\hrulefill}\thinspace}
\providecommand{\MR}{\relax\ifhmode\unskip\space\fi MR }
% \MRhref is called by the amsart/book/proc definition of \MR.
\providecommand{\MRhref}[2]{%
  \href{http://www.ams.org/mathscinet-getitem?mr=#1}{#2}
}
\providecommand{\href}[2]{#2}

\end{document}